\documentclass[11pt]{article}
\usepackage{tikz}
\usepackage{graphicx}
\usepackage{psfrag}
\usepackage{mathrsfs}
\usepackage{amsfonts}
\usepackage{amsmath}
\usepackage[center]{caption2}
\usepackage{amssymb}
\usepackage{fancyhdr}
\usepackage{color}
\usepackage[top=1.54cm,bottom=2.54cm,left=2.54cm,right=2.54cm]{geometry}
\usepackage{indentfirst}
\def\v{\varepsilon}

\def\q{\quad}
\def\qq{\quad\quad\quad\quad\quad}

\def \B{\begin{equation}}
\def\bbb{\begin{array}{lll}}

\def\eee{\end{array}}
\def \E{\end{equation}}

\def\it{\indent}

\begin{document}
\title{Temporal periodic solutions to nonhomogeneous quasilinear hyperbolic equations driven by time-periodic boundary conditions}
\author{
{\small  Xixi Fang$^{1}$, Peng Qu $^{2}$, Huimin Yu$^{3}$,
        \footnote{Email: fangxixi\_a@163.com}
         \footnote{Email: pqu@fudan.edu.cn.}
          \footnote{Corresponding author, Email: hmyu@sdnu.edu.cn}}\\
 {\small \it $^{1,3}$ School of Mathematics and Statistics, Shandong Normal University, Jinan 250014, China,}\\
 {\small \it $^{2}$ School of Mathematical Sciences and Shanghai Key Laboratory for Contemporary Applied Mathematics,}\\
 {\small \it Fudan University, Shanghai 200433, China}}
\date{}
\maketitle

\begin{center}
\begin{minipage}{160mm}
{\bf Abstract}\hskip 3mm \it We consider the temporal periodic solutions to general nonhomogeneous quasilinear hyperbolic equations with
 a kind of weak diagonal dominant structure. Under the temporal periodic boundary conditions, the existence, stability and uniqueness of the time-periodic classical solutions are obtained. Moreover, the $W^{2,\infty}$ regularity and stability around the time-periodic solutions are discussed.
 Our results reveal that the feedback boundary control with dissipative structure can stabilize the $K$-weakly diagonally dominant nonhomogeneous quasilinear hyperbolic system around the temporal periodic solution.
\end{minipage}
\end{center}
\section{Introduction}
\setcounter{equation}{0}
\it In this paper, we consider the initial-boundary value problem of a kind of nonhomogeneous quasilinear hyperbolic equations
\begin{equation}\label{1.1}
\partial_t u+A(u) \partial_x u= F(u), \quad(t, x) \in \mathbb{R} \times[0, L],
\end{equation}
where $u=\left(u_1(t, x), \ldots, u_n(t, x)\right)^T \in C^1(\mathbb{R} \times[0, L]; \mathcal{U})$ is the unknown vector,
$\mathcal{U} \subset \mathbb{R}^n $ is a small neighborhood of $u=0$,
$n\times n$ matrix $A(u)=\left(a_{i j}(u)\right)_{i, j=1}^n$ and
$n$-dimensional vector $F(u)=\left(f_1(u), \ldots, f_n(u)\right)^T$ are both defined on $\mathcal{U}$.
The assumptions on the coefficient matrix $A(u)=\left(a_{i j}(u)\right)_{i, j=1}^n$  are the same as \cite{Qu}.
For the completeness of this paper, we reiterate these symbols here.
The smooth matrix-valued function $A(u)$ has $n$ nonzero real eigenvalues $\lambda_i(u)(i=1, \ldots, n)$
and the corresponding smooth left and right eigenvectors
\begin{equation}\label{1.2}
l_i(u)=(l_{i1}(u),...,l_{in}(u)), \q r_i(u)=(r_{i1}(u),...,r_{in}(u))^{T}, \q \forall u \in \mathcal{U}, \forall i=1, \ldots, n
\end{equation}
satisfying
\begin{equation}\label{1.3}
\lambda_r(u)<0<\lambda_s(u), \quad  \forall u \in \mathcal{U}, \forall r=1, \ldots, m ; s=m+1, \ldots, n,
\end{equation}
and
\begin{equation}\label{1.4}
l_i(u) A(u)=\lambda_i(u) l_i(u),
\q\q
 A(u) r_i(u)=\lambda_i(u) r_i(u), \quad  \forall u \in \mathcal{U}, \forall i=1, \ldots. n.
\end{equation}
Without loss of generality, we may assume that
\begin{equation}\label{1.6}
A(0)=\operatorname{diag}\left\{\lambda_i(0)\right\}_{i=1}^n
\end{equation}
is a diagonal matrix.
Moreover, by changing the length of $l_i(u),r_j(u)$, we continue to assume that

\begin{equation}\label{1.7}
l_i(u) r_j(u)=\delta_{i j}, \q  \forall u \in \mathcal{U}, \forall i, j=1, \ldots, n,
\end{equation}
\begin{equation}\label{1.8}
\left|r_i(u)\right|=1,      \q\q  \forall u \in \mathcal{U}, \forall i=1, \ldots, n
\end{equation}
for simplicity.
Here $\delta_{i j}$ is the Kronecker's symbol, and thus
\begin{equation}\label{1.9}
l_i(0)=e_i^T, \quad r_i(0)=e_i, \quad l_{i j}(0)=r_{i j}(0)=\delta_{i j},  \q \forall i, j=1, \ldots, n
\end{equation}
if we denote the $i$-th unit vector in $\mathbb{R}^n$ by $e_i$.
Notice $\lambda_i (u)$ $(i=1,...,n)$ are nonzero, we set
\begin{equation}\label{1.11}
\mu_i(u)=\lambda_i^{-1}(u),
\end{equation}
and
\begin{equation}\label{1.12}
\mu_{\max }=\max _{i=1, \ldots, n} \sup _{u \in \mathcal{U} }\left|\mu_i(u)\right| .
\end{equation}
By rescaling the time variable (if needed), we can further assume
\begin{equation}\label{1.13}
\mu_{\max } \leq 1 .
\end{equation}
We assume the nonhomogeneous term $F(u)$ satisfies
\begin{equation}\label{1.14}
F(0)=0
\end{equation}
and $\nabla F(0)$ satisfies the \emph{K-weakly diagonally dominant} condition,
that is, for
$$\mathcal{G}(u)=\nabla F(u)$$
or
\begin{equation*}
g_{ij}(u)=\frac{\partial}{\partial u_j}f_i(u),\q \forall i,j=1,...,n,
\end{equation*}
there exists  a positive constant $K\geq0$ such that
\begin{equation}\label{1.15}
(-g_{ii}(0))-(\sum_{j\neq i}|g_{ij}(0)|)>-K.
\end{equation}
\\
\textbf{Remark 1.}\emph{
If system $(\ref{1.1})$ is homogeneous (i.e. $F(u)=0$), or the nonhomogeneous term is superlinear (i.e. $\nabla F(0)=0$), then $(\ref{1.15})$ is satisfied  obviously.}\\\\
\textbf{Remark 2.}\emph{
If system $(\ref{1.1})$ meets the strictly diagonally dominant condition defined in Hsiao\cite{Hsiao}, i.e.,
$$
(-g_{ii}(0)) - (\sum_{j\neq i}|g_{ij}(0)|) \geq b >0,
$$
for a positive constant $b$. Then $(\ref{1.15})$ holds and $K$ is taken as zero.}\\

We consider ($\ref{1.1}$) with initial data
\begin{equation}\label{1.16'}
u(x,0) =u_0(x), \q\q  x\in [0, L],
\end{equation}
and boundary conditions
\begin{align}
& x=0: u_s=G_s\left(h_s(t), u_1, \ldots, u_m\right), \q\q\q   s=m+1, \ldots, n,  \label{1.16}\\
& x=L: u_r=G_r\left(h_r(t), u_{m+1}, \ldots, u_n\right), \q\q   r=1, \ldots, m.\q   \label{1.17}
\end{align}
Here  we assume $G_s\left(h_s, u_1, \ldots, u_m\right)(s=m+1, \ldots, n)$ and $G_r\left(h_r, u_{m+1}, \ldots, u_n\right)(r=1, \ldots, m)$ are smooth functions with
\begin{align}
& G_s(0,0, \ldots, 0)=0, \quad s=m+1, \ldots, n,\label{1.18} \\
& G_r(0,0, \ldots, 0)=0, \quad r=1, \ldots, m,\label{1.19}
\end{align}
and $h_i(t)$ ($i=1,...,n$) are small $C^1$ smooth functions.
We further suppose that
\begin{equation}\label{1.20}
\max _{r=1, \ldots, m}\left|\frac{\partial G_r}{\partial h_r}(0, \ldots, 0)\right| \leq \frac{1}{2}, \q
\max _{s=m+1, \ldots, n}\left|\frac{\partial G_s}{\partial h_s}(0, \ldots, 0)\right| \leq \frac{1}{2},
\end{equation}
otherwise, if
\begin{equation*}
\max _{r=1, \ldots, m}\left|\frac{\partial G_r}{\partial h_r}(0, \ldots, 0)\right| \leq M_0, \q \text{and} \q M_0 > \frac{1}{2},
\end{equation*}
define
\begin{equation*}
\tilde{h}_r =2 M_0 h_r,
\end{equation*}
then
\begin{equation*}
\max _{r=1, \ldots, m}\left|\frac{\partial G_r}{\partial \tilde{h}_r}(0, \ldots, 0)\right| \leq \frac{1}{2}.
\end{equation*}

As \cite{Qu}, we also need the boundary conditions ($\ref{1.16}$)-($\ref{1.17}$) to be dissipative, that is, for the matrix
$$
\Theta=\left(\theta_{i j}\right)_{i, j=1}^n \stackrel{\text { def. }}{=}
\left(\begin{array}{cc}
0 & \left(\frac{\partial G_r}{\partial u_s}(0,0, \ldots, 0)\right)_{\substack{r=1, \ldots, m \\
s=m+1, \ldots, n}} \\
\left(\frac{\partial G_s}{\partial u_r}(0,0, \ldots, 0)\right)_{\substack{s=m+1, \ldots, n \\
r=1, \ldots, m}} & 0
\end{array}\right),
$$
we need its minimal characterizing number satisfies
\begin{equation}\label{1.23}
\theta=\|\Theta\|_{\min }=\inf _{\substack{\Gamma=\operatorname{diag}\left\{\gamma_i\right\}_{i=1}^n \\ \gamma_i \neq 0}}\left\|\Gamma \Theta \Gamma^{-1}\right\|<1,
\end{equation}
where
$$
\|\Theta\|=\max _{i=1, \ldots, n} \sum_{j=1}^n\left|\theta_{i j}\right|.
$$

In the past decades, many efforts have been made for homogeneous or inhomogeneous hyperbolic conservation laws, such as \cite{Hsiao,Liu,Qin,YingWang}.
In 1982, Hsiao \cite{Hsiao} promoted a simple algorithm for the nonhomogeneous hyperbolic systems of balance laws with strictly diagonally dominant dissipation.
Kawashima and Shizuta \cite{Kawa,Shi-Kawa} put forward the famous Kawashima condition, which had been widely applied in the following studies by many experts. The authors showed that the strict dissipativity brings about the decay of solutions.
After that, many related results emerged, we can see \cite{Kawa1,Kawa-Shi,Kawa-Yong,Plaza,Shi-Kawa,Xu,Yong} and the references therein.
On the other hand, the dissipative boundary conditions are needed for considering classical solutions to the initial-boundary problem of quasilinear hyperbolic system, see \cite{Li} for details.
In 1985, Qin \cite{Qin} considered the homogeneous quasilinear hyperbolic systems with dissipative boundary conditions.
Coron and Bastin \cite{Coron} studied the exponential stability of classical solutions of one-dimensional quasi-linear hyperbolic systems with dissipative boundary conditions. The authors showed the Lyapunov stability results for $C^1$ norm.
At the same time, the authors investigated the $BV$ stability of $2\times2$ hyperbolic systems of conservation laws with strictly positive velocities in  \cite{Coron1}.
In 2015, Li and Liu \cite{LL} considered the combined effect of the internal dissipative condition and the boundary dissipative condition for inhomogeneous quasilinear hyperbolic systems with small initial data.

As for time-periodic solutions of partial differential equations, there are many contributions on viscous fluid equations \cite{CT,Luo,MUY,MN} and hyperbolic conservation laws \cite{Nao,OS,Takeno}.
However, all the time-periodic solutions they concerned about are induced by time-periodic external forces.
As for the time-periodic solutions derived from time-periodic boundary conditions of hyperbolic conservation law, there remains few results.
Yuan\cite{Yuan} made a breakthrough in considering the existence and high-frequency limiting behavior of supersonic time-periodic solutions for one-dimensional isentropic compressible Euler flows. It is worthy to be point out that the temporal periodic solutions considered in \cite{Yuan} are driven by periodic boundary conditions.
For a more general quasilinear hyperbolic systems, Qu \cite{Qu} proved the existence, uniqueness and stability of time-periodic solutions arising from the dissipative and time-periodic boundary.
Then, for the one-dimensional isentropic compressible Euler equations with a linear damping, \cite{Qu1} considered the time-periodic solutions triggered by time-periodic boundary in a bounded domain.
\cite{Zhang2} learned the isentropic compressible Euler equation with nonlinear damping.
As for non isentropic compressible Euler equations, Temple and Young \cite{TY,TY1} proved the existence of space and time periodic solutions, we can see them for more details.

In this paper, we are interested in the periodic solutions derived by the time-periodic boundary conditions with a common constant temporal period $T_*>0$,
\begin{equation}\label{1.22}
h_i\left(t+T_*\right)=h_i(t), \quad \forall t \in \mathbb{R},\q \forall i=1, \ldots, n.
\end{equation}
We aim to extend all the results in \cite{Qu} for homogeneous quasilinear hyperbolic equations to general nonhomogeneous systems with $K$-weakly diagonally condition $(\ref{1.15})$.

The main results are as follows.
\\\\
\textbf{Theorem 1} \emph{(Existence of time-periodic solutions).
Under the assumptions $(\ref{1.15})$ and $(\ref{1.23})$, for any small enough $K>0$, there exists a small constant $\varepsilon_1 > 0$ and a constant $C_E >0$, such that for any given $\varepsilon \in\left(0, \varepsilon_1\right)$, any given $T_* \in \mathbb{R}_{+}$, and any given functions $h_i(t)(i=1, \ldots, n)$ satisfying $(\ref{1.22})$ and
\begin{equation}\label{1.24}
\max_{i=1,...,n} \|h_i\|_{C^1} \leq \varepsilon,
\end{equation}
then there exists a $C^1$ smooth function $u_0(x)$ with
\begin{equation}\label{1.24'}
\|u_0\|_{C^1} \leq C_E\varepsilon,
\end{equation}
such that system $(\ref{1.1})$ and $(\ref{1.16'})$-$(\ref{1.17})$ admits a $C^1$ classical solution $u=$ $u^{(P)}(t, x)$ in the domain
$D=\{(t, x)| t\in \mathbb{R}, x \in[0, L]\}$, satisfying
\begin{equation}\label{1.25}
u^{(P)}\left(t+T_*, x\right)=u^{(P)}(t, x), \quad \forall(t, x) \in D
\end{equation}
and
\begin{equation}\label{1.26}
\|u^{(P)}\|_{C^1(D)} \leq C_E\varepsilon.
\end{equation}}
\\
\textbf{Remark 3.}\emph{
In fact, from the proof of Theorem 1, the constant $K$ does not have to be arbitrary small, we can determine an upper bound for $K$. The upper bound depends only on $\theta$, $L$ and the maximum of weight functions $W_i(x)$ $(i=1,...,n)$. However, this would lead to a tedious elementary calculation which is not the main purpose of this paper.
}\\
\\
\textbf{Theorem 2 }\emph{(Stability of the time-periodic solution).
Under the assumptions $(\ref{1.15})$ and $(\ref{1.23})$, for any arbitrary small $K>0$, there exists a small constant $\varepsilon_2 \in\left(0, \varepsilon_1\right)$ and a constant $C_S>0$, such that for any given $\varepsilon \in\left(0, \varepsilon_2\right)$ and any given functions $u_0(x)$ and $h_i(t)(i=1, \ldots, n)$ satisfying $(\ref{1.22})$-$(\ref{1.24'})$ with certain compatibilities,
system $(\ref{1.1})$ and $(\ref{1.16'})$-$(\ref{1.17})$ admits a unique global $C^1$ classical solution $u=u(t, x)$ on the domain $D_1=\{(t, x)| t\in \mathbb{R}_+, x \in[0, L]\}$,
satisfying
\begin{equation}\label{1.27}
\left\|u(t, \cdot)-u^{(P)}(t, \cdot)\right\|_{C^0} \leq C_S \varepsilon \beta^{[t/T_0]}, \quad \forall t \in \mathbb{R}_{+},
\end{equation}
where $u^{(P)}$, depending on $h_i(t)(i=1, \ldots, n)$, is the time-periodic solution obtained in Theorem 1, $\beta\in(0,1)$  and
\begin{equation}\label{1.28}
T_0= \max_{i=1,...,n} \sup_{u \in \mathcal{U} } \frac{L}{\left|\lambda_i(u)\right|}=L \mu_{max} .
\end{equation}}
\\
Let $t\rightarrow +\infty$ in $(\ref{1.27})$, then we have the uniqueness result.\\
\\
\textbf{Corollary 3} \emph{(Uniqueness of the time-periodic solution).
Under the assumptions $(\ref{1.15})$ and $(\ref{1.23})$, there exists a small constant $\varepsilon_3 \in (0, \varepsilon_2)$, such that for any given $\varepsilon \in (0, \varepsilon_3)$ and any given functions $h_i(t)(i=1, \ldots, n)$ satisfying $(\ref{1.22})$ and $(\ref{1.24})$, the time-periodic solution $u=$ $u^{(P)}(t, x)$ obtained in Theorem 1 is unique.}
\\\\
\textbf{Theorem 4 }\emph{(Regularity of the time-periodic solution).
Under the assumptions $(\ref{1.15})$ and $(\ref{1.23})$,  for all functions $h_i(t)$ satisfying $(\ref{1.22})$ and  $(\ref{1.24})$ with further $W^{2,\infty}$ regularity as
\begin{equation}\label{1.29}
\max_{i=1,...,n} \|h''_i\|_{L^\infty} \leq M_0 < +\infty,
\end{equation}
there exist constants $\varepsilon_4 \in\left(0, \varepsilon_1\right)$ and $M_R>0$, such that for any given $\varepsilon \in\left(0, \varepsilon_4\right)$,
the time-periodic solution $u=u^{(P)}(t, x)$  obtained in Theorem 1 on the domain $D_1=\{(t, x)| t\in \mathbb{R}_+, x \in[0, L]\}$,
also obeys the $W^{2,\infty}$ regularity with
\begin{equation}\label{1.30}
\max_{i=1,...,n}\left\{ \|\partial_t^2 u^{(P)}\|_{L^\infty}, \|\partial_t \partial_x u^{(P)}\|_{L^\infty},
\|\partial_x^2 u^{(P)}\|_{L^\infty},\right\} \leq M_R < +\infty.
\end{equation}
}
\\
Based on Theorem 4, we can further get the stability around the time-periodic solution, that is,\\
\\
\textbf{Theorem 5 }\emph{(Exponential stability around regular time-periodic solution).
Under the same assumptions in Theorem 2 and further $W^{2,\infty}$ regularity  $(\ref{1.29})$,
there exist constants $\varepsilon_5 \in\left(0, \min\{\varepsilon_2, \varepsilon_4\}\right)$ and $\widetilde{C}_S>0$,
such that for any given $\varepsilon \in\left(0, \varepsilon_5\right)$, we have further stability around the time-periodic solution
in $C^1$ with \begin{equation}\label{1.31}
\max\{\left\|\partial_t u(t, \cdot)-\partial_t u^{(P)}(t, \cdot)\right\|_{C^0}, \left\|\partial_x u(t, \cdot)-\partial_x u^{(P)}(t, \cdot)\right\|_{C^0}\} \leq \widetilde{C}_S \varepsilon \beta^{[t/T_0]}, \quad \forall t \in \mathbb{R}_{+}.
\end{equation}
}

Based on a ingeniously constructed linear iterative scheme and meticulous estimates for the iterative sequence, we will prove Theorems 1-5 in Sections 2-5.
\section{Existence of time-periodic solutions}
\setcounter{equation}{0}

In this section, we will construct a time-periodic solutions to problem ($\ref{1.1}$) and ($\ref{1.16'}$)-($\ref{1.17}$) by an iteration method. To this end, we apply the classical linearization as in \cite{Hu,LY}.  Multiply system ($\ref{1.1}$) by the left eigenvectors $l_i (u)$ ($i=1,...,n$) from left to get
\begin{equation}\label{2.1}
\partial_t u_i+\lambda_i(u) \partial_x u_i-f_i(u)=\sum_{j=1}^n B_{i j}(u)\left(\partial_t u_j+\lambda_i(u) \partial_x u_j-f_j(u)\right),\q i=1,...,n, \\
\end{equation}
where
\begin{equation}\label{2.2}
\begin{array}{ll}
B_{i j}(u)=
\begin{cases}
-\displaystyle\frac{l_{i j}(u)}{l_{i i}(u)}, & i \neq j \\
\q 0, & i=j
\end{cases}
\end{array}
\end{equation}
satisfying
\begin{equation}\label{2.3}
B_{ij}(0)=0, \quad \forall i, j=1,...n.
\end{equation}
To change the roles of $x$ and $t$, we multiply $\mu_i(u)= \lambda_i^{-1}(u)$ to the $i-th$ equation of ($\ref{2.1}$) and get
\begin{equation}\label{2.4}
\begin{aligned}
   \partial_x u_i+\mu_i(u) \partial_t u_i- & \mu_i(u) f_i(u)\\
 =& \sum_{j=1}^n B_{i j}(u)\left(\partial_x u_j+\mu_i(u)\partial_t u_j-\mu_i(u)f_j(u)\right), \quad \forall i, j=1,...n.
\end{aligned}
\end{equation}
Define
\begin{equation}\label{2.6}
\begin{array}{ll}
\tilde{g}_{i j}=
\begin{cases}
 \q\q\mu_i(0)\displaystyle\frac{\partial f_i}{\partial u_j}(0)&=\mu_i(0) g_{i j}(0), \q\q\q\q  j\neq i, \\
 \mu_i(0)\left(\displaystyle\frac{\partial f_i}{\partial u_i}(0)-K \right)&=\mu_i(0)\left(g_{i i}(0)-K\right), \quad j=i,
\end{cases}
\end{array}
\end{equation}
and
\begin{equation}\label{2.7}
 \tilde{g}_i^{N L}\left(u\right)
=\mu_i(u) f_i(u)-\sum_{j=1}^n \mu_i(0) \displaystyle\frac{\partial f_i}{\partial u_j}(0) u_j-\sum_{j=1}^n B_{i j}(u)\mu_i(u) f_j(u).
\end{equation}
Then it is easy to get
\begin{equation}\label{2.8}
\tilde{g}_i^{N L}(0)=0, \q \nabla_u \tilde{g}_i^{N L}(0)=0,
\end{equation}
i.e., $\tilde{g}_i^{N L}$ is a nonlinear term
and ($\ref{2.4}$) can be rewritten into
\begin{equation}\label{2.8'}
\begin{aligned}
\partial_x u_i+\mu_i\left(u\right) \partial_t u_i-\tilde{g}_{i i} u_i
& =\sum_{j=1}^n B_{i j}\left(u\right)\left(\partial_x u_j+\mu_i\left(u\right) \partial_t u_j\right) \\ & +\sum_{j\neq i} \tilde{g}_{i j} u_j+K \mu_i(0) u_i+\tilde{g}_i^{N L}\left(u\right).
\end{aligned}
\end{equation}
Utilize the $K$-weakly diagonally dominant type condition ($\ref{1.15}$), and notice the fact
$$\quad-1 \leq-\mu_{\max } \leq \mu_r(0)<0, \quad \forall r=1, \ldots, m,$$
we have
$$-\tilde{g}_{r r}+\left(\sum_{j\neq r}\left|\tilde{g}_{r j}\right|\right)=\left|\mu_r(0)\right| \cdot\left(g_{r r}(0)-K\right)+\sum_{j\neq r}\left|\mu_r(0)\right|\left|g_{r j}(0)\right| < 0,
\q \forall r=1, \ldots, m. $$
So we get $\tilde{g}_{rr}$ satisfies the diagonally dominate condition
\begin{equation}\label{2.9}
\tilde{g}_{r r} >\big(\sum_{j \neq r}\left|\tilde{g}_{r j}\right|\big) \geq 0, \q \forall r=1,...,m.
\end{equation}
Similarly, since
$$0<\mu_{s(0)} \leq \mu_{\max } \leq 1,\q \forall s=m+1,...,n,$$
we have
\begin{equation}\label{2.10}
\quad-\tilde{g}_{s s}>\big(\sum_{j \neq s}\left|\tilde{g}_{s j}\right|\big) \geq 0 , \quad \forall s=m+1,...,n.
\end{equation}

Now we consider the following iterative system
\begin{equation}\label{2.5}
\begin{aligned}
\partial_x u_i^{(l)}+\mu_i\left(u^{(l-1)}\right) \partial_t u_i^{(l)}-\tilde{g}_{i i} u_i^{(l)}
& =\sum_{j=1}^n B_{i j}\left(u^{(l-1)}\right)\left(\partial_x u_j^{(l-1)}+\mu_i\left(u^{(l-1)}\right) \partial_t u_j^{(l-1)}\right) \\ & +\sum_{j\neq i} \tilde{g}_{i j} u_j^{(l-1)}+K \mu_i(0) u_i^{(l-1)}+\tilde{g}_i^{N L}\left(u^{(l-1)}\right)
\end{aligned}
\end{equation}
for $l \in \mathbb{Z}_{+} $, and start our iteration from
\begin{equation}\label{2.5'}
u^{(0)} (t,x) \equiv 0
\end{equation}
with the  linearized boundary conditions
\begin{equation}\label{2.11}
\q\    x=0: u_s^{(l)}=G_s\left(h_s(t), u_1^{(l-1)}, \ldots, u_m^{(l-1)}\right), \quad s=m+1, \ldots, n,
\end{equation}
\begin{equation}\label{2.12}
x=L: u_r^{(l)}=G_r\left(h_r(t), u_{m+1}^{(l-1)}, \ldots, u_n^{(l-1)}\right), \quad r=1, \ldots, m.
\end{equation}
We first prove the following priori estimates :
\\\\
\textbf{Proposition 2.1.} \emph{Suppose the boundary hypotheses satisfy $(\ref{1.20})$-$(\ref{1.24})$, for small enough $K>0$, the sequence of $C^{1}$ solutions $u_i^{(l)}(t,x)$ $(i=1,...,n; l \in \mathbb{Z}_{+})$ to the linearized iteration system $(\ref{2.5})$ and $(\ref{2.11})$-$(\ref{2.12})$ starting from $(\ref{2.5'})$ satisfy
\begin{flalign}\label{2.17}
\begin{split}
\it \left.1 \right) \q\q\q
u_i^{(l)}\left(t+T_\ast, x\right)=u_i^{(l)}(t, x),\q \forall (t, x)\in \mathbb{R}\times[0,L],
\end{split}&
\end{flalign}
\begin{equation}\label{2.13}
\left\|u_i^{(l)}\right\|_{C^{1}} \leq M_1 \varepsilon,
\end{equation}
and
\begin{equation}\label{2.14}
\left\|u_i^{(l)}-u_i^{(l-1)}\right\|_{C^0} \leq M_2 \beta^l \varepsilon
\end{equation}
for $\beta\in(0,1)$, small enough $\varepsilon >0$ and large enough constants $M_1, M_2>0$.\\
\it $\left.2 \right)$
Moreover, we have
\begin{equation}\label{2.15}
\max_{i=1,...,n}\left\{2 \omega\left(\delta \mid \partial_t u_i^{(l)}(\cdot, x)\right), \omega\left(\delta \mid \partial_x u_i^{(l)}(\cdot, x)\right)\right\} \leq \frac{1}{4} \Omega(\delta),
\end{equation}
\begin{equation}\label{2.16}
\max_{i=1,...,n}\left\{\omega\left(\delta \mid \partial_t u_i^{(l)}\right), \omega\left(\delta \mid \partial_x u_i^{(l)}\right)\right\} \leq \Omega(\delta),
\end{equation}
where $\omega(\delta|\cdot)$ denote modulus of continuity in some sense, by definitions of
\begin{equation*}
\omega(\delta \mid f(. , x))=\sup _{\substack{\left|t_1-t_2\right| \leq \delta \\}}\left|f\left(t_1, x\right)-f\left(t_2, x\right)\right|
\end{equation*}
and
\begin{equation*}
\omega(\delta \mid f)=\sup _{\substack{\left|t_1-t_2\right| \leq \delta \\\left|x_1-x_2\right| \leq \delta}}\left|f\left(t_1, x_1\right)-f\left(t_2, x_2\right)\right|,
\end{equation*}
$\Omega(\delta)$ denotes a continuous function of $\delta \in (0,1)$, independent of $l$, and to be defined later in $(\ref{omega})$
with
$$\lim_{\delta\rightarrow 0^+} \Omega(\delta)=0.$$
}

Resemble \cite{Qu},  the proof of Theorem 1 can be derived by Proposition 2.1 directly.
First, from $(\ref{2.14})$, we can get the convergence of the sequence $\{u^{(l)}\}_{l=1}^{\infty}$
according to the $C^0$ Cauchy sequence property, and thus converges to a $C^0$ function $u^{(P)}$ uniformly.
Then by $(\ref{2.17})$, we show that $u^{(P)}$ is time-periodic.
After that, with $(\ref{2.13})$ and $(\ref{2.16})$, we can utilize the Arzel\`{a}-Ascoli theorem to get the existence of a subsequence of $\{u_i^{(l)}\}$ that converges to $u^{(P)}$ in $C^{1}$ uniformly.
By uniqueness, $\{u_i^{(l)}\}$ $(i=1,...,n)$ converges to $u^{(P)}$ in $C^{1}$ norm.
Therefore, the limit of sequence, $u^{(P)}$, is a $C^{1}$ smooth solution to system $(\ref{1.1})$ and $(\ref{1.16'})$-$(\ref{1.17})$.
To get the initial data, we take $u^{(P)}(x) = u^{(P)}(0,x)$, then from $(\ref{2.13})$ we can get $(\ref{1.24'})$.
\\\\
\textbf{The proof of Proposition 2.1 :} We will prove $(\ref{2.17})$-$(\ref{2.16})$ in Proposition 2.1 inductively.  Firstly, it is easy to see that $(\ref{2.17})$-$(\ref{2.13})$  and $(\ref{2.15})$-$(\ref{2.16})$ are satisfied naturally for $u_i^{(0)}=0$,
then for each $l\in \mathbb{Z}_{+}$,
$i=1,...,n,$
we give the following assumptions
\begin{equation}\label{2.22}
u_i^{(l-1)}\left(t+T_\ast, x\right)=u_i^{(l-1
)}(t, x), \q \forall (t, x)\in \mathbb{R}\times[0,L],
\end{equation}
\begin{equation}\label{2.18}
\left\|u_i^{(l-1)}\right\|_{C^{1}} \leq M_1 \varepsilon,
\end{equation}
\begin{equation}\label{2.19}
\left\|u_i^{(l-1)}-u_i^{(l-2)}\right\|_{C^0} \leq M_2 \beta^{l-1} \varepsilon,\q ( \text{for }  l\geq 2)
\end{equation}
and
\begin{equation}\label{2.20}
\max_{i=1,...,n}\left\{2 \omega\left(\delta \mid \partial_t u_i^{(l-1)}(\cdot, x)\right), \omega\left(\delta \mid \partial_x u_i^{(l-1)}(\cdot, x)\right)\right\} \leq \frac{1}{4} \Omega(\delta),
\q \forall x\in[0,L],
\end{equation}
\begin{equation}\label{2.21}
\max_{i=1,...,n}\left\{\omega\left(\delta \mid \partial_t u_i^{(l-1)}\right), \omega\left(\delta \mid \partial_x u_i^{(l-1)}\right)\right\} \leq \Omega(\delta).
\end{equation}

Since ($\ref{2.5}$) and ($\ref{2.11}$)-($\ref{2.12}$) are decoupled nonhomogeneous linear transport equations for $u_i^{(l)}$ $(i=1,...,n, l\in \mathbb{Z}_{+})$,  and notice the time-periodic boundary condition $(\ref{1.22})$ and our initial value of iteration ($\ref{2.5'}$), the periodic conclusion ($\ref{2.17}$) can be obtained by the uniqueness of solution.

Now, we prove the $C^0$ estimates in ($\ref{2.13}$). From the boundary conditions ($\ref{2.11}$), for any $s=m+1,...,n$, note
($\ref{1.20}$) and ($\ref{1.23}$), we have
\begin{equation}\label{2.23}
\begin{aligned}
\left|u_s^{(l)}(t, 0)\right| & =\left|\int_0^1\left(\frac{\partial G_s}{\partial h_s}\left(\tau \cdot h_s, \tau \cdot u_r\right) h_s(t)\right.\right.\\
 & \q\q  \left.\left.+\sum_{r=1}^m \frac{\partial G_s}{\partial u_r}\left(\tau \cdot h_s, \tau \cdot u_r\right) u_r^{(l-1)}(t, 0)\right) d \tau\right| \\
& \leq \frac{\varepsilon}{2}+\theta M_1 \varepsilon .
\end{aligned}
\end{equation}
Similarly, on $x=L$, we have
\begin{equation}\label{2.25}
\left|u_r^{(l)}(t, L)\right| \leq \frac{\varepsilon}{2}+\theta M_1 \varepsilon,
\q\q \forall r=1,...,m .
\end{equation}
In order to get the $C^0$ estimates in the domain $\{(t,x)| t \in \mathbb{R} ,0<x< L\}$, we define the characteristic curve $t=t_i^{(l)}(x;t_0,x_0)$ for $i=1,...,n$
and $l\in \mathbb{Z}_{+}$ by the following ODEs:
\begin{equation}\label{2.26}
\left\{\begin{array}{l}
\displaystyle\frac{d}{d x} t_i^{(l)}\left(x ; t_0, x_0\right)=\mu_i\left(u^{(l-1)}\left(t_i^{(l)}\left(x ; t_0, x_0\right), x\right)\right), \\\\
t_i^{(l)}\left(x_0 ; t_0, x_0\right)=t_0.
\end{array}\right.
\end{equation}
Define weight functions $W_i(x)$ $(i=1,...,n)$ as
\begin{equation}\label{2.27}
W_r(x)=e^{\tilde{g}_{r r}(L-x)}, \q W_s(x)=e^{-\tilde{g}_{s s} x} ,
\end{equation}
we can easily get
\begin{align}
 W_r(x) \geq 1, \q W_s(x) \geq 1 ,& \q\q W_r(L)=1, \q W_s(0)=1,\label{2.28*}\\
  \nonumber  \\
 W_r^{\prime}(x)=-\tilde{g}_{r r} W_r(x) <0, & \q\q W_s^{\prime}(x)=-\tilde{g}_{s s} W_s(x)>0,\label{2.28}
\end{align}
thus we have
\begin{align}
 1 \leq W_r(x)\leq \max\{W_r(0), W_s(L)\}=M_3, \q\q
 1 \leq W_s(x)\leq M_3, \label{M3}
\end{align}
and
\begin{equation}\label{2.29}
\int_L^{x_0} W_r(x) d x=\frac{1-W_r\left(x_0\right)}{\tilde{g}_{r r}} \leq 0, \q
\int_0^{x_0} W_s(x) d x=\frac{1-W_s\left(x_0\right)}{\tilde{g}_{s s}} \geq 0 .
\end{equation}
Note definitions of $W_r (x)$ and $W_s(x)$ in $(\ref{2.27})$,
we transform equation ($\ref{2.5}$) into
\begin{equation}\label{2.30}
\begin{aligned}
  &\partial_x \left( W_r u_r^{(l)}\right)+\mu_r\left(u^{(l-1)}\right) \partial_t \left(W_r u_r^{(l)}\right)\\
= & W_r\left(\partial_x u_r^{(l)}+\mu_r\left(u^{(l-1)}\right) \partial_t u_r^{(l)}-\tilde{g}_{r r } u_r^{(l)}\right)\\
= & \sum_{j=1}^n W_r \cdot B_{r j}\left(u^{(l-1)}\right)\left(\partial_x u_j^{(l-1)}+\mu_r\left(u^{(l-1)}\right)
    \partial_t u_j^{(l-1)}\right) + W_r \widetilde{g}_r^{N L}(u^{(l-1)})\\
  & +K W_r \mu_r(0) u_r^{(l-1)}+\sum_{j\neq r} W_r\tilde{g}_{r j} u_j^{(l-1)},\q
  \q\q r=1,...,m
\end{aligned}
\end{equation}
and
\begin{equation}\label{2.31}
\begin{aligned}
  &\partial_x \left( W_s u_s^{(l)}\right)+\mu_s\left(u^{(l-1)}\right) \partial_t \left(W_s u_s^{(l)}\right)\\
= & W_s\left(\partial_x u_s^{(l)}+\mu_s\left(u^{(l-1)}\right) \partial_t u_s^{(l)}-\tilde{g}_{s s} u_s^{(l)}\right)\\
= & \sum_{j=1}^n W_s \cdot B_{s j}\left(u^{(l-1)}\right)\left(\partial_x u_j^{(l-1)}+\mu_s\left(u^{(l-1)}\right)
    \partial_t u_j^{(l-1)}\right) + W_s \widetilde{g}_s^{N L}(u^{(l-1)})\\
  & +K W_s \mu_s(0) u_s^{(l-1)}+\sum_{j\neq s} W_s\tilde{g}_{s j} u_j^{(l-1)},\q
  \q\q s=m+1,...,n.
\end{aligned}
\end{equation}
For $s=m+1,...,n$, we integrate ($\ref{2.31}$) along the characteristic curve $t=t_s^{(l)}(x;t_0,x_0)$ from $x=0$ to $x=x_0$
to get
\begin{equation}\label{2.32}
\begin{aligned}
  & \ W_s(x_0) u_s^{(l)}(t_0, x_0)-W_s(0) u_s^{(l)}(t_s^{(l)}(0;t_0,x_0), 0)\\
= & \int_0^{x_0}\left(\sum_{j=1}^n W_s B_{s j}\left(u^{(l-1)}\right)\left(\partial_x u_j^{(l-1)}+\mu_s\left(u^{(l-1)}\right)
    \partial_t u_j^{(l-1)}\right)+W_s \tilde{g}_s^{N L}(u^{(l-1)}) \right.\\
    &\left. +K W_s \mu_s(0) u_s^{(l-1)}+\sum_{j\neq s} W_s\tilde{g}_{s j} u_j^{(l-1)}\right) dx.
\end{aligned}
\end{equation}
Noticing that $\sum\limits_{j=1}^n W_s B_{s j}\left(u^{(l-1)}\right)\left(\partial_x u_j^{(l-1)}+\mu_s\left(u^{(l-1)}\right)
    \partial_t u_j^{(l-1)}\right)$ and $W_s \tilde{g}_s^{N L}(u^{(l-1)})$ are nonlinear terms, from ($\ref{2.32}$) we have
\begin{equation}\label{2.33}
\begin{aligned}
\left|u_s^{(l)}(t_0, x_0)\right|
& \leq \frac{ \frac{\varepsilon}{2} + \theta M_1 \varepsilon}{W_s\left(x_0\right)}
   +\frac{C\varepsilon^2 + K L W_s(L) M_1 \varepsilon }{W_s\left(x_0\right)}
   +\frac{ \left(\sum\limits_{j\neq s}|\tilde{g}_{sj}|\right)\left(W_s(x_0)-1\right)} {|\tilde{g}_{ss}| W_s\left(x_0\right)} \cdot M_1\varepsilon\\
& \leq \frac{1}{W_s\left(x_0\right)}
  \bigg( \varepsilon+\theta M_1 \varepsilon + K L M_3\cdot M_1 \varepsilon - M_1 \varepsilon  \bigg)
  + M_1 \varepsilon,
\end{aligned}
\end{equation}
where we have used the condition $(\ref{2.10})$ in the second inequality.\\
If
\begin{equation}\label{2.33*}
\theta + K L M_3 <1, \q\text{and}\q   M_1 > \frac{1}{1-\theta-K L M_3},
\end{equation}
then ($\ref{2.33}$) implies
\begin{equation}\label{2.33**}
\left|u_s^{(l)}(t_0, x_0)\right| \leq  M_1 \varepsilon.
\end{equation}
Here
$$\theta + K L M_3 <1$$
denote the coupled relationship between dissipation boundary and $K$-weakly diagonally dominant source term.
Similarly, for $r=1,...,m$, integrate ($\ref{2.30}$) along the characteristic curve $t=t_r^{(l)}(x;t_0,x_0)$ from $x=L$ to $x=x_0$ to get
\begin{equation}\label{2.34'}
\begin{aligned}
  & \ W_r(x_0) u_r^{(l)}(t_0, x_0)-W_r(L) u_r^{(l)}(t_r^{(l)}(L;t_0,x_0), L)\\
= & \int_L^{x_0}\left(\sum_{j=1}^n W_r B_{r j}\left(u^{(l-1)}\right)\left(\partial_x u_j^{(l-1)}+\mu_r\left(u^{(l-1)}\right)
    \partial_t u_j^{(l-1)}\right)
    + W_r \tilde{g}_r^{N L}(u^{(l-1)}) \right.\\
    &\left. +K W_r \mu_r(0) u_r^{(l-1)}+\sum_{j\neq r} W_r\tilde{g}_{r j} u_j^{(l-1)}\right) dx,
\end{aligned}
\end{equation}
then from $(\ref{2.9})$, $(\ref{2.25})$ and $(\ref{2.33*})$, we get
\begin{equation}\label{2.34}
\begin{aligned}
\left|u_r^{(l)}(t_0, x_0)\right|
& \leq \frac{ \frac{\varepsilon}{2} + \theta M_1 \varepsilon}{W_r\left(x_0\right)}
   +\frac{C\varepsilon^2 + K L W_r(0) M_1 \varepsilon }{W_r\left(x_0\right)}
   +\frac{\left(\sum\limits_{j\neq r}|\tilde{g}_{rj}|\right) \left(W_r(x_0)-1\right) } {|\tilde{g}_{rr}| W_r\left(x_0\right)} \cdot M_1\varepsilon\\
& \leq \frac{1}{W_r\left(x_0\right)}
  \bigg( \varepsilon+\theta M_1 \varepsilon +K L M_3\cdot M_1 \varepsilon - M_1 \varepsilon  \bigg)
  + M_1 \varepsilon \\
& \leq M_1 \varepsilon.\\
\end{aligned}
\end{equation}
Then we finish the $C^0$ norm estimates.
As noted in Remark 3, for convenience, we will always assume $K$ sufficiently small in the following proof.

For $C^1$ estimates, we denote
\begin{equation}\label{2.35}
z_i^{(l)}=\partial_t u_i^{(l)}, \q \omega_i^{(l)}=\partial_x u_i^{(l)},\q i=1,...,n,
\end{equation}
then on the boundary, applying $\partial_t$  to ($\ref{2.11}$) and ($\ref{2.12}$) yields
\begin{equation}\label{2.36}
\begin{aligned}
\quad x=0: z_s^{(l)}
&=\frac{\partial G_s}{\partial h_s}\left(h_s, u_1^{(l-1)},..., u_m^{(l-1)}\right) h_s^{\prime}(t)\\
&+\sum_{r=1}^m \frac{\partial G_s}{\partial u_r}\left(h_s, u_1^{(l-1)}, \ldots, u_m^{(l-1)}\right) z_r^{(l-1)},
\q  \forall s=m+1,...,n,
\end{aligned}\end{equation}
\begin{equation}\label{2.37}
\begin{aligned}
x=L: z_r^{(l)}
&=\frac{\partial G_r}{\partial h_r}\left(h_r, u_{m+1}^{(l-1)},\ldots, u_n^{(l-1)}\right) h_r^{\prime}(t)\\
&+\sum_{s=m+1}^n \frac{\partial G_r}{\partial u_s}\left(h_s, u_{m+1}^{(l-1)}, \ldots, u_r^{(l-1)}\right) z_s^{(l-1)},
\q  \forall r=1,...,m.
\end{aligned}\end{equation}
From ($\ref{1.20}$), ($\ref{1.23}$) and ($\ref{1.24}$) we deduce
\begin{equation}\label{2.38}
\left|z_s^{(l)}(t, 0)\right| \leq \frac{\varepsilon}{2}+\theta M_1 \varepsilon \leq \alpha M_1 \varepsilon,
\q \forall s=m+1,...,n,
\end{equation}
and
\begin{equation}\label{2.39}
\left|z_r^{(l)}(t, 0)\right| \leq \frac{\varepsilon}{2}+\theta M_1 \varepsilon \leq \alpha M_1 \varepsilon,
\q  \forall r=1,...,m,
\end{equation}
where we take $\alpha=\frac{\theta+1}{2}\in (\theta,1)$ and $M_1 >\frac{100}{1-\alpha}> 0$.
In the domain  $\{(t,x)| t \in \mathbb{R} ,0<x< L\}$, we take the temporal derivative to ($\ref{2.5}$) and get
\begin{equation}\label{2.40}
\begin{aligned}
  & \partial_x z_i^{(l)}+\mu_i\left(u^{(l-1)}\right) \partial_t z_i^{(l)}-\tilde{g}_{i i} z_i^{(l)}\\
= & -\sum_{j=1} ^n  \displaystyle\frac{\partial \mu_i}{\partial u_j}\left(u^{(l-1)}\right) z_j^{(l-1)} z_i^{(l)}
  +\sum_{j, k=1} ^n  B_{i j}\left(u^{(l-1)}\right) \frac{\partial \mu_i}{\partial u_k}\left(u^{(l-1)}\right) z_k^{(l-1)} z_j^{(l-1)}\\
  & +\sum_{j, k=1}^n \frac{\partial B_{i j}}{\partial u_k}\left(u^{(l-1)}\right) z_k^{(l-1)}\left(\omega_j^{(l-1)}+\mu_i\left(u^{(l-1)}\right) z_j^{(l-1)}\right)\\
  & +\sum_{j=1}^n \left(\partial_x+\mu_i\left(u^{(l-1)}\right) \partial_t\right)\left(B_{i j}\left(u^{(l-1)}\right) z_j^{(l-1)}\right)\\
  & -\sum_{j , k=1}^n  \frac{\partial B_{i j}}{\partial u_k}\left(u^{(l-1)}\right) z_j^{(l-1)}\left(\omega_k^{(l-1)}+\mu_i\left(u^{(l-1)}\right) z_k^{(l-1)}\right) \\
  & +\sum_{j=1}^n \frac{\partial \tilde{g}_i^{N L}}{\partial u_j}\left(u^{(l-1)}\right) z_j^{(l-1)}
  +K \mu_i(0) z_i^{(l-1)}+\sum_{j \neq i} \tilde{g}_{i j} z_j^{(l-1)},
\q \forall i=1,...,n.
\end{aligned}\end{equation}
Moreover, we rewrite ($\ref{2.40}$) with $W_r(x)$ $(r=1,...,m)$ and $W_s(x)$ $(s=m+1,...,n)$ in the form
$$\begin{aligned}
   & \quad \partial_x\left(W_r z_r^{(l)}\right)+\mu_r\left(u^{(l-1)}\right) \partial_t\left(W_r z_r^{(l)}\right) \qq\qq\qq\\
= & -\sum_{j=1}^n \frac{\partial \mu_r}{\partial u_j}\left(u^{(l-1)}\right) z_j^{(l-1)}\left(W_r z_r^{(l)}\right)
\end{aligned}$$
\begin{equation}\label{2.41}
\begin{aligned}
  & +\sum_{j=1}^n \left(\partial_x+\mu_r\left(u^{(l-1)}\right) \partial_t\right)\left(W_r B_{r j}\left(u^{(l-1)}\right) z_j^{(l-1)}\right)\\
  & +W_r \left( \sum_{j , k=1}^n  \frac{\partial B_{r j}}{\partial u_k}\left(u^{(l-1)}\right)\left(z_k^{(l-1)} \omega_j^{(l-1)}-z_j^{(l-1)} \omega_k^{(l-1)}\right)\right.\\
  &  +\sum_{j, k=1}^n  B_{r j}\left(u^{(l-1)}\right) \frac{\partial \mu_r}{\partial u_k}\left(u^{(l-1)}\right) z_k^{(l-1)} z_j^{(l-1)}
     +\sum_{j=1}^n \tilde{g}_{r r} B_{r j}\left(u^{(l-1)}\right) z_j^{(l-1)}\\
  & \left. +\sum_{j=1}^n \frac{\partial \tilde{g}_r^{N L}}{\partial u_j}\left(u^{(l-1)}\right) z_j^{(l-1)}
   +K \mu_r(0) z_r^{(l-1)}
  +\sum_{j \neq r} \tilde{g}_{r j} z_j^{(l-1)}\right),
\end{aligned}\end{equation}
and
\begin{equation}\label{2.42}
\begin{aligned}
   & \quad \partial_x\left(W_s z_s^{(l)}\right)+\mu_s\left(u^{(l-1)}\right) \partial_t\left(W_s z_s^{(l)}\right) \\
= & -\sum_{j=1}^n \frac{\partial \mu_s}{\partial u_j}\left(u^{(l-1)}\right) z_j^{(l-1)}\left(W_s z_s^{(l)}\right) \\
  & + \sum_{j=1}^n \left(\partial_x+\mu_s\left(u^{(l-1)}\right) \partial_t\right)\left(W_s B_{s j}\left(u^{(l-1)}\right) z_j^{(l-1)}\right) \\
  & +W_s \left( \sum_{j , k=1}^n \frac{\partial B_{s j}}{\partial u_k}\left(u^{(l-1)}\right)\left(z_k^{(l-1)} \omega_j^{(l-1)}-z_j^{(l-1)} \omega_k^{(l-1)}\right)\right.\\
  & +\sum_{j, k=1}^n B_{s j}\left(u^{(l-1)}\right) \frac{\partial \mu_s}{\partial u_k}\left(u^{(l-1)}\right) z_k^{(l-1)} z_j^{(l-1)}
    + \sum_{j=1}^n \tilde{g}_{s s} B_{s j}\left(u^{(l-1)}\right) z_j^{(l-1)}\\
    & \left. +\sum_{j=1}^n \frac{\partial \tilde{g}_s^{N L}}{\partial u_j}\left(u^{(l-1)}\right) z_j^{(l-1)}
      +K \mu_s(0) z_s^{(l-1)}
     +\sum_{j \neq s} \tilde{g}_{s j} z_j^{(l-1)}\right).
\end{aligned}\end{equation}
Notice that $\nabla \tilde{g}_i^{N L}(0)=0$, the term $\sum\limits_{j=1}^n \frac{\partial \tilde{g}_i^{N L}}{\partial u_j}\left(u^{(l-1)}\right) z_j^{(l-1)}$  is equivalent to $\varepsilon^2$.
To get the $C^0$ estimates of $z_r^{(l)}$ and $z_s^{(l)}$ by ($\ref{2.41}$) and ($\ref{2.42}$), the main difficulty comes from the first term of the right hand. We overcome this aporia by \emph{bootstrap arguments}.

For $s=m+1,...,n$, exchange the role of $x$ and $t$, $z_s^{(l)}$ satisfy the decoupled linear transport equation ($\ref{2.36}$)($\ref{2.42}$) if we regard $u^{(l-1)}$ and $z^{(l-1)}$ as given functions.
Therefore, for any $(t,x)\in \mathbb{R}\times[0,L]$, we have the existence and uniqueness of $z_s^{(l)}(t,x)$.
Moreover, note the initial condition $(\ref{2.36})$, there exists a constant $L^* \in(0, L]$ such that
\begin{equation}\label{2.43}
\max _{x \in\left[0, L^*\right]} \max_{t}\left|z_s^{(l)}(t, x)\right| \leq M_1 \varepsilon.
\end{equation}
Now we show that $L^*$ can be extended to $L$ by contradiction.
If ($\ref{2.43}$) holds, and
\begin{equation}\label{2.44}
\max _t\left|z_s^{(l)}\left(t, L^*\right)\right|=M_1 \varepsilon,
\end{equation}
we will present that
\begin{equation}\label{2.44'}
\max _{x \in\left[0, L^*\right]} \max _t\left|z_s^{(l)}(t, x)\right|<M_1 \varepsilon
\end{equation}
for any $L^{*}\leq L$.
Integrate ($\ref{2.42}$) along the characteristic curve $t=t_s^{(l)}(x;t_0,x_0)$ from $x=0$ to $x=x_0\leq L^{*}$,
and we get\\
\begin{equation}\label{2.45}
\begin{aligned}
  & \q W_s\left(x_0\right) z_s^{(l)}\left(t_0, x_0\right)-W_s(0) z_s^{(l)}\left(t_s^{(l)}\left(0 ; t_0, x_0\right), 0\right) \\
 =& \int_0^{x_0}\left(-\left.\sum_{j=1}^n \frac{\partial \mu_s}{\partial u_j}\left(u^{(l-1)}\right) z_j^{(l-1)} \left(W_s z_s^{(l)}\right)\right)\right|_{\left(t_s^{(l)}\left(x_; t_0, x_0\right), x\right)} d x \\
  & +\left.\left(W_s B_{s j}\left(u^{(l-1)}\right) z_j^{(l-1)}\right)\right|_{\left(t_{s}^{l} \left(0; t_0, x_0\right), 0 \right)} ^{\left(t_0, x_0\right)} \\
  & +\int_0^{x_0} W_s \left.
  \left( \sum_{j , k=1}^n \frac{\partial B_{s j}}{\partial u_k}\left(u^{(l-1)}\right)
  \left(z_k^{(l-1)} \omega_j^{(l-1)}-z_j^{(l-1)} \omega_k^{(l-1)}\right)\right.\right.\\
  & +\sum_{j, k=1}^n B_{s j}\left(u^{(l-1)}\right) \frac{\partial \mu_s}{\partial u_k}\left(u^{(l-1)}\right) z_k^{(l-1)} z_j^{(l-1)}
  +\sum_{j=1}^n \tilde{g}_{s s} B_{s j}\left(u^{(l-1)}\right) z_j^{(l-1)}\\
  & \left.\left.  +\sum_{j=1}^n \frac{\partial \tilde{g}_s^{N L}}{\partial u_j}\left(u^{(l-1)}\right) z_j^{(l-1)}+K \mu_s(0) z_s^{(l-1)}\right)
  \right|_{\left(t_s^{(l)}\left(x_; t_0, x_0\right), x\right)}  d x\\
  &+\left.\int_0^{x_0}\left(W_s \sum_{j \neq s} \tilde{g}_{s j} z_j^{(l-1)}\right)\right|_{\left(t_s^{l}\left(x; t_0, x_0\right), x\right)} d x, \\
  &
\end{aligned}
\end{equation}
then we have
\begin{equation}\label{2.46}
\begin{aligned}
\left|z_s^{(l)}\left(t_0, x_0\right)\right|
\leq \ & \frac{\alpha M_1 \varepsilon }{W_s\left(x_0\right)}
            +\frac{C\left(\varepsilon^2 + K \varepsilon\right)}{W_s\left(x_0\right)}
            +\frac{ \bigg(\sum\limits_{j\neq s}|\tilde{g}_{sj}|\bigg)  \left(W_s(x_0)-1\right)} {|\tilde{g}_{ss}| W_s\left(x_0\right)} \cdot M_1\varepsilon\\
\leq \ & \frac{\left(M_1-100\right)\varepsilon}{W_s\left(x_0\right)}
       +\frac{\varepsilon}{W_s\left(x_0\right)}
       +\frac{W_s\left(x_0\right)-1}{W_s\left(x_0\right)} \cdot M_1 \varepsilon \\
 =   \ &  M_1 \varepsilon-\frac{99}{W_s\left(x_0\right)} \varepsilon
 \ <   \   M_1 \varepsilon,
\end{aligned}
\end{equation}
here we have used the fact $M_1 >\frac{100}{1-\alpha}> 0$.
Consequently, for any $(t,x)\in \mathbb{R}\times[0,L]$, we get the temporal derivative estimation
\begin{equation}\label{2.47}
\left|z_s^{(l)}\left(t, x\right)\right|\leq M_1 \varepsilon,
\q\q \forall s=m+1,...,n.
\end{equation}
Similarly, integrating $(\ref{2.41})$ along the characteristic curve $t=t_r^{(l)}(x;t_0,x_0)$ from $x=L$ to $x=x_0\geq 0$, we obtain
\begin{equation}\label{2.48}
\left|z_r^{(l)}\left(t, x\right)\right|\leq M_1 \varepsilon,
\q\q \forall r=1,...,m, \forall (t,x)\in \mathbb{R}\times[0,L].
\end{equation}
Furthermore, take $M_1 > 2 \max\limits_{i=1,...,n} \sum_{j=1}^n \left|\tilde{g}_{i j}\right|$, then from ($\ref{2.5}$) we can directly derive
 \begin{equation}\label{2.49}
\left|\omega_i^{(l)}\left(t, x\right)\right|\leq M_1 \varepsilon ,
\q\q \forall i=1,...,n.
\end{equation}
Therefor, the $C^1$ estimates are obtained.

Next, we will prove the $C^0$ Cauchy sequence property ($\ref{2.14}$).
For $l=1$, since ($\ref{2.5'}$), it can be derived from ($\ref{2.33}$)-($\ref{2.34}$) directly.
For $l\geq 2$, on the boundary, we have
\begin{equation}\label{2.50}
\begin{aligned}
\q\q\q\q x=0: \ \ & u_s^{(l)}-u_s^{(l-1)} \\
       = & \sum_{r=1}^m \int_0^1 \frac{\partial G_s}{\partial u_r}\left(h_s(t),
                                  \tau u_1^{(l-1)}+(1-\tau) u_1^{(l-2)},\ldots,\right.\\
          & \left. \tau u_m^{(l-1)}+(1-\tau) u_m^{(l-2)}\right) d \tau
            \cdot\left(u_r^{(l-1)}-u_r^{(l-2)}\right),
            \q \forall s=m+1,...,n,
\end{aligned}
\end{equation}
\begin{equation}\label{2.51}
\begin{aligned}
x=L: \q& u_r^{(l)}-u_r^{(l-1)} \\
   = & \sum_{s=m+1}^n \int_0^1 \frac{\partial G_r}{\partial u_s}\left(h_r(t), \tau u_{m+1}^{(l-1)}+(1-\tau) u_{m+1}^{(l-2)}, \ldots,\right.\\
    & \left.\tau u_n^{(l-1)}+(1-\tau) u_n^{(l-2)}\right) d \tau
     \cdot\left(u_s^{(l-1)}-u_s^{(l-2)}\right),
     \q \forall r=1,...,m
\end{aligned}
\end{equation}
from ($\ref{2.11}$)-($\ref{2.12}$), and then
\begin{equation}\label{2.52}
\left|u_s^{(l)}(t, 0)-u_s^{(l-1)}(t, 0)\right| \leq \theta\left(1+C \varepsilon \right) M_2 \beta^{l-1} \varepsilon ,
\end{equation}
\begin{equation}\label{2.53}
\left|u_r^{(l)}(t, L)-u_r^{(l-1)}(t, L)\right| \leq \theta\left(1+C \varepsilon \right) M_2 \beta^{l-1} \varepsilon .
\end{equation}
Notice that $W_i$ $(i=1,...,n)$ are independent on $l$, then from ($\ref{2.30}$)-($\ref{2.31}$) we arrive at \\
\begin{equation}\label{2.54}
\begin{aligned}
  &\q \partial_x\left(W_r\left(u_r^{(l)}-u_r^{(l-1)}\right)\right)+\mu_r\left(u^{(l-1)}\right) \partial_t\left(W_r\left(u_r^{(l)}-u_r^{(l-1)}\right)\right) \\
= & -W_r \partial_t u_r^{(l-1)}\left(\mu_r\left(u^{(l-1)}\right)-\mu_r\left(u^{(l-2)}\right)\right) \quad \\
  & +W_r\left(\sum_{j=1}^n B_{r j}\left(u^{(l-1)}\right)\left(\partial_x+\mu_r\left(u^{(l-1)}\right) \partial_t\right)\left(u_j^{(l-1)}-u_j^{(l-2)}\right)\right.\\
  & +\sum_{j=1}^n\left(B_{r j}\left(u^{(l-1)}\right)-B_{r j}\left(u^{(l-2)}\right)\right)\left(\partial_x+\mu_r\left(u^{(l-1)}\right) \partial_t\right) u_j^{(l-2)} \\
  & +\sum_{j=1}^n B_{r j}\left(u^{(l-2)}\right)\left(\mu_r\left(u^{(l-1)}\right)-\mu_r\left(u^{(l-2)}\right)\right) \partial_t u_j^{(l-2)} +\tilde{g}_r^{N L}\left(u^{(l-1)}\right)\\
  & -\tilde{g}_r^{N L}\left(u^{(l-2)}\right) +K \mu_r(0)\left(u_r^{(l-1)}-u_r^{(l-2)}\right)
   \left.+\sum_{j\neq r} \tilde{g}_{r j}\left(u_j^{(l-1)}-u_j^{(l-2)}\right)\right)
\end{aligned}
\end{equation}
$$
\begin{aligned}
= & \left(\partial_x+\mu_r\left(u^{(l-1)}\right) \partial_t\right)\left(\sum_{j=1}^n W_r B_{r j}\left(u^{(l-1)}\right)\left(u_j^{(l-1)}-u_j^{(l-2)}\right)\right) + I_1 \\
  & + W_r \left(\beta \sum_{j \neq r} \tilde{g}_{r j}\left(u_j^{(l-1)}-u_j^{(l-2)}\right)+(1-\beta) \sum_{j \neq r} \tilde{g}_{r j}\left(u_j^{(l-1)}-u_j^{(l-2)}\right)\right) \\
\end{aligned}
$$
and
\begin{equation}\label{2.55}
\begin{aligned}
  &\q \partial_x\left(W_s\left(u_s^{(l)}-u_s^{(l-1)}\right)\right)+\mu_s\left(u^{(l-1)}\right) \partial_t\left(W_s\left(u_s^{(l)}-u_s^{(l-1)}\right)\right) \\
= & \left(\partial_x+\mu_s\left(u^{(l-1)}\right) \partial_t\right)\left(\sum_{j=1}^n W_s B_{s j}\left(u^{(l-1)}\right)\left(u_j^{(l-1)}-u_j^{(l-2)}\right)\right) +I_2 \\
  & + W_s \left(\beta \sum_{j \neq s} \tilde{g}_{s j}\left(u_j^{(l-1)}-u_j^{(l-2)}\right)+(1-\beta) \sum_{j \neq s} \tilde{g}_{s j}\left(u_j^{(l-1)}-u_j^{(l-2)}\right)\right), \\
\end{aligned}
\end{equation}
where
$$
\begin{aligned}
I_1 = & W_r\bigg(-\sum_{j=1}^n \partial_t u_r^{(l-1)} \int_0^1 \frac{\partial \mu_r}{\partial u_j}\left(\tau u^{(l-1)}+(1-\tau) u^{(l-2)}\right) d \tau \cdot\left(u_j^{(l-1)}-u_j^{(l-2)}\right) \\
  & -\sum_{j,k=1}^n \frac{\partial B_{r j}}{\partial u_k}\left(u^{(l-1)}\right)\left(\partial_x u_k^{(l-1)}+\mu_r\left(u^{(l-1)}\right) \partial_t u_k^{(l-1)}\right)\left(u_j^{(l-1)}-u_j^{(l-2)}\right) \\
  & +\sum_{j=1}^{n} B_{r j}\left(u^{(l-1)}\right) \tilde{g}_{r r}\left(u_j^{(l-1)}-u_j^{(l-2)}\right)
   +\sum_{j,k=1}^{n}\left(\partial_x u_j^{(l-2)} +\mu_r\left(u^{(l-1)}\right) \partial_t u_j^{(l-2)}\right)\\
   & \q\q\q\q \cdot\int_0^1 \frac{\partial B_{r j}}{\partial u_k}\left(\tau u^{(l-1)}+(1-\tau) u^{(l-2)}\right) d \tau \cdot\left(u_k^{(l-1)}-u_k^{(l-2)}\right)
\end{aligned}
$$
$$
\begin{aligned}
  & +\sum_{j, k=1}^n  B_{r j}\left(u^{(l-1)}\right) \partial_t u_j^{(l-2)} \int_0^1 \frac{\partial \mu_r}{\partial u_k}\left(\tau u^{(l-1)}+(1-\tau) u^{(l-2)}\right) d \tau \cdot\left(u_k^{(l-1)}-u_k^{(l-2)}\right)\\
  & +\sum_{j=1}^n  \int_0^1 \frac{\partial \tilde{g}_r^{N L}}{\partial u_j}\left(\tau u_j^{(l-1)}+(1-\tau) u_j^{(l-2)}\right) d \tau \cdot\left(u_j^{(l-1)}-u_j^{(l-2)}\right) \\
  &  +K \mu_r(0)\left(u_r^{(l-1)}-u_r^{(l-2)}\right) \bigg),\\
\end{aligned}
$$
$$\begin{aligned}
I_2 =
& W_s \bigg( -\sum_{j=1}^n \partial_t u_s^{(l-1)} \int_0^1 \frac{\partial \mu_s}{\partial u_j}\left(\tau u^{(l-1)}+(1-\tau) u^{(l-2)}\right) d \tau \cdot\left(u_j^{(l-1)}-u_j^{(l-2)}\right) \\
  & -\sum_{j,k=1}^n \frac{\partial B_{s j}}{\partial u_k}\left(u^{(l-1)}\right)\left(\partial_x u_k^{(l-1)}+\mu_s\left(u^{(l-1)}\right) \partial_t u_k^{(l-1)}\right)\left(u_j^{(l-1)}-u_j^{(l-2)}\right) \\
  & +\sum_{j=1}^{n} B_{s j}\left(u^{(l-1)}\right) \tilde{g}_{ss}\left(u_j^{(l-1)}-u_j^{(l-2)}\right)
    +\sum_{j,k=1}^{n}\left(\partial_x u_j^{(l-2)}+\mu_s\left(u^{(l-1)}\right) \partial_t u_j^{(l-2)}\right)\\
    &\q\q\q\q  \cdot\int_0^1 \frac{\partial B_{s j}}{\partial u_k}\left(\tau u^{(l-1)}+(1-\tau) u^{(l-2)}\right) d \tau \cdot\left(u_k^{(l-1)}-u_k^{(l-2)}\right) \\
  & +\sum_{j, k=1}^n  B_{s j}\left(u^{(l-1)}\right) \partial_t u_j^{(l-2)} \int_0^1 \frac{\partial \mu_s}{\partial u_k}\left(\tau u^{(l-1)}+(1-\tau) u^{(l-2)}\right) d \tau \cdot\left(u_k^{(l-1)}-u_k^{(l-2)}\right) \\
  & +\sum_{j=1}^n  \int_0^1 \frac{\partial \tilde{g}_s^{N L}}{\partial u_j}\left(\tau u_j^{(l-1)}+(1-\tau) u_j^{(l-2)}\right) d \tau \cdot\left(u_j^{(l-1)}-u_j^{(l-2)}\right) \\
  &  +K \mu_s(0)\left(u_s^{(l-1)}-u_s^{(l-2)}\right)\bigg). \\
\end{aligned}$$

Integrate ($\ref{2.55}$) along the characteristic curve $t=t_s^{(l)}(x;t_0,x_0)$ $(s=m+1,...,n)$ from $x=0$ to $x=x_0$ to get
\begin{equation}\label{2.56}
\begin{aligned}
  & W_s\left(x_0\right)\left(u_s^{(l)}\left(t_0, x_0\right)-u_s^{(l-1)}\left(t_0, x_0\right)\right)\\
  & -W_s(0)\left(u_s^{(l)}\left(t_s^{(l)}\left(0; t_0, x_0\right), 0\right)
  -u_s^{(l-1)}\left(t_s^{(l)}\left(0 ; t_0, x_0\right), 0\right)\right) \\
= & \left.\sum_{j=1}^n\left( W_s B_{s j}\left(u^{(l-1)}\right)\left(u_j^{(l-1)}-u_j^{(l-2)}\right)\right)
    \right|
    _{\left(t_s^{\left(l\right)}\left(0; t_0, x_0\right), 0 \right)}
    ^{(t_0, x_0)}
   +\int_0^{x_0} I_2 \left.\right|_{t=t_s^{(l)}\left(x ; t_0, x_0\right)} d x \\
  & +\left.\beta \int_0^{x_0} \sum_{j \neq s} W_s \tilde{g}_{sj}\left(u_j^{(l-1)}-u_j^{(l-2)}\right)\right|_{t=t_s^{(l)}\left(x_; t_0, x_0\right)} d x \\
  & +(1-\beta) \int_0^{x_0} \sum_{j \neq s} \left. W_s \tilde{g}_{s j} \left(u_j^{(l-1)}-u_j^{(l-2)}\right)\right|_{t=t_s^{(l)}\left(x_; t_0, x_0\right)} d x, \\
\end{aligned}
\end{equation}
thus we have
\begin{equation}\label{2.57}
\begin{aligned}
   & \left|u_s^{(l)}\left(t_0, x_0\right)-u_s^{(l-1)}\left(t_0, x_0\right)\right|\\
 \leq & \frac{\theta (1+ C \varepsilon)}{W_s (x_0)} M_2 \beta^{l-1}\varepsilon
        +\frac{C (\varepsilon + K)}{W_s (x_0)} M_2 \beta^{l-1}\varepsilon
        +\frac{\left( \sum\limits_{j\neq s}\left|\tilde{g}_{sj}\right| \right)\left( W_s(x_0)-1 \right)}{\left| \tilde{g}_{ss} \right|W_s (x_0)} M_2 \beta^{l}\varepsilon \\
      & +\frac{\left( \sum\limits_{j\neq s}\left|\tilde{g}_{sj}\right| \right)\left( W_s(x_0)-1 \right)}{\left| \tilde{g}_{ss} \right|W_s (x_0)} \left(1-\beta\right)M_2 \beta^{l-1}\varepsilon \\
 \leq & M_2 \beta^l \varepsilon+\frac{M_2 \beta^{l-1} \varepsilon}{W_s\left(x_0\right)}
      \bigg(\theta\left(1+C \varepsilon \right)
      +C\left(\varepsilon+K\right)+\left(1-\beta\right)W_s\left(x_0\right)-1 \bigg).
\end{aligned}
\end{equation}
Take $0<1-\beta \ll 1, \beta \in(0,1)$ such that
\begin{equation}\label{beta}
1-\beta<\min _{i=1,...,n} \inf_{x} \frac{1- \theta}{W_i(x)},
\end{equation}
where $\beta$ is independent of $t, x, l$ and $\varepsilon$.
Take $\varepsilon \ll 1$, $K \ll 1$, then we have
$$\theta(1+C \varepsilon)+C(\varepsilon+K )+(1-\beta) W_s\left(x_0\right)-1<0.$$
Therefore,
\begin{equation}\label{2.58}
\left|u_s^{(l)}\left(t_0, x_0\right)-u_s^{(l-1)}\left(t_0, x_0\right)\right| \leq M_2 \beta^{l}\varepsilon,
\q \forall s=m+1,...,n.
\end{equation}
Similarly we have
\begin{equation}\label{2.59}
\left|u_r^{(l)}\left(t_0, x_0\right)-u_r^{(l-1)}\left(t_0, x_0\right)\right| \leq M_2 \beta^{l}\varepsilon,
\q \forall r=1,...,m.
\end{equation}
Thus, the proof of ($\ref{2.14}$) is finished.

Now we focus on  ($\ref{2.15}$) and  ($\ref{2.16}$). First we show the modulus of continuity for $u_i^{(l)}$ on the temporal directions, i.e. ($\ref{2.15}$).
On the boundary $x=0$, for any two points $(t_1, 0)$ and $(t_2, 0)$ with $|t_1 -t_2| \leq \delta$, from ($\ref{2.36}$) we have
\begin{equation}\label{2.60}
\begin{aligned}
   &  z_s^{(l)}\left(t_2, 0\right)-z_s^{(l)}\left(t_1, 0\right) \\
 = & \frac{\partial G_s}{\partial h_s}\left(h_s, u_1^{(l-1)}, \ldots, u_m^{(l-1)}\right)\left(t_2, 0\right) \left(h_s^{\prime}\left(t_2\right)-h_s^{\prime}\left(t_1\right)\right) \\
& +\int_0^1\left( \frac { \partial^{2} G_{s} } { \partial h_{ s }^{ 2 } }
           \left(\tau h_s\left(t_2\right)+(1-\tau) h_s\left(t_1\right),
                 \tau u_1^{(l-1)}\left(t_2, 0\right)+(1-\tau) u_1^{(l-1)}\left(t_1, 0\right),\ldots,
           \right.\right.\\
&  \left. \q
   \tau u_m^{(l-1)}\left(t_2, 0\right)+(1-\tau) u_m^{(l-1)}\left(t_1, 0\right)\right)
   \cdot \left(h_s\left(t_2\right)-h_s\left(t_1\right)\right) \\
&  +\sum_{r=1}^m \frac{\partial^2 G_s}{\partial u_r \partial h_s}
  \left(\tau h_s\left(t_2\right)+(1-\tau) h_s\left(t_1\right),
  \tau u_1^{(l-1)}\left(t_2, 0\right)+(1-\tau) u_1^{(l-1)}(t_1, 0), \cdots,\right. \\
&  \q \tau u_m^{(l-1)}\left(t_2, 0\right)+(1-\tau) u_m^{(l-1)}\left(t_1, 0\right) \Big)
   \cdot\left(u_r^{(l-1)}\left(t_2, 0\right)-u_r^{(l-1)}\left(t_1, 0\right)\right)\Bigg)
   d \tau \cdot h_s^{'}\left(t_1\right) \\
& +\sum_{r=1}^m \frac{\partial G_s}{\partial u_r}
   \left(h_s, u_1^{(l-1)}, \ldots, u_m^{(l-1)}\right)\left(t_2, 0\right)
   \cdot\left(z_r^{(l-1)}\left(t_2, 0\right)-z_r^{(l-1)}\left(t_1, 0\right)\right)\\
& +\int_0^1\left(\frac{\partial^2 G_s}{\partial h_s \partial u_r}
           \left(\tau h_s\left(t_2\right)+(1-\tau) h_s\left(t_1\right),
           \tau u_1^{(l-1)}\left(t_2, 0\right)+(1-\tau) u_1^{(l-1)}\left(t_1, 0\right), \ldots,
           \right.\right.\\
&  \left.\q
   \tau u_m^{(l-1)}\left(t_2, 0\right)+(1-\tau) u_m^{(l-1)}\left(t_1, 0\right)\right)
   \cdot  \left(h_s\left(t_2\right)-h_s(t_1)\right) \\
& +\sum_{\tilde{r}=1}^m \frac{\partial^2 G_s}{\partial u_{\tilde{r}} \partial u_r}
  \left(\tau h_s\left(t_2\right)+(1-\tau) h_s\left(t_1\right),
  \tau u_1^{(l-1)}\left(t_2, 0\right)+(1-\tau) u_1^{(l-1)}(t_1, 0), \cdots,\right. \\
&  \q \tau u_m^{(l-1)}\left(t_2, 0\right)+(1-\tau) u_m^{(l-1)}\left(t_1, 0\right) \Big)
   \cdot\left(u_{\tilde{r}}^{(l-1)}\left(t_2, 0\right)-u_{\tilde{r}}^{(l-1)}\left(t_1, 0\right)\right)\Bigg)
   d \tau \\
&  \q \cdot z_r^{(l-1)}\left(t_1, 0\right),
\end{aligned}
\end{equation}
 and for any two points $(t_2, L)$ and $(t_1, L)$ with $|t_1 -t_2| \leq \delta$ on $x=L$, we have
\begin{equation}\label{2.61}
\begin{aligned}
 &  z_r^{(l)}\left(t_2, L\right)-z_r^{(l)}\left(t_1, L\right) \\
 = & \frac{\partial G_r}{\partial h_r}\left(h_r, u_{m+1}^{(l-1)},
     \ldots, u_n^{(l-1)}\right)\left(t_2, L\right) \left(h_r^{\prime}\left(t_2\right)-h_r^{\prime}\left(t_1\right)\right) \\
& +\int_0^1\left( \frac { \partial^{2} G_{r} } { \partial h_{ r }^{ 2 } }
           \left(\tau h_r\left(t_2\right)+(1-\tau) h_r\left(t_1\right),
                 \tau u_{m+1}^{(l-1)}\left(t_2, L\right)+(1-\tau) u_{m+1}^{(l-1)}\left(t_1, L\right),\ldots,
           \right.\right.\\
&  \left. \q
   \tau u_n^{(l-1)}\left(t_2, L\right)+(1-\tau) u_n^{(l-1)}\left(t_1, L\right)\right)
   \cdot \left(h_r\left(t_2\right)-h_r\left(t_1\right)\right) \\
&  +\sum_{s=m+1}^n \frac{\partial^2 G_r}{\partial u_s \partial h_r}
  \left(\tau h_r\left(t_2\right)+(1-\tau) h_r\left(t_1\right),
  \tau u_{m+1}^{(l-1)}\left(t_2, L\right)+(1-\tau) u_{m+1}^{(l-1)}(t_1, L), \cdots,\right. \\
&  \q \tau u_n^{(l-1)}\left(t_2, L\right)+(1-\tau) u_n^{(l-1)}\left(t_1, L\right) \Big)
   \cdot\left(u_r^{(l-1)}\left(t_2, L\right)-u_r^{(l-1)}\left(t_1, L\right)\right)\bigg)
   d \tau \cdot h_r^{'}\left(t_1\right) \\
& +\sum_{s=m+1}^n \frac{\partial G_r}{\partial u_s}
   \left(h_r, u_{m+1}^{(l-1)}, \ldots, u_n^{(l-1)}\right)\left(t_2, L\right)
   \cdot\left(z_r^{(l-1)}\left(t_2, L\right)-z_r^{(l-1)}\left(t_1, L\right)\right)
\end{aligned}
\end{equation}
$$\begin{aligned}
& +\int_0^1\left(\frac{\partial^2 G_r}{\partial h_r \partial u_s}
           \left(\tau h_r\left(t_2\right)+(1-\tau) h_r\left(t_1\right),
           \tau u_{m+1}^{(l-1)}\left(t_2, L\right)+(1-\tau) u_{m+1}^{(l-1)}\left(t_1, L\right), \ldots,
           \right.\right.\\
&  \left.\q
   \tau u_n^{(l-1)}\left(t_2, L\right)+(1-\tau) u_n^{(l-1)}\left(t_1, L\right)\right)
   \cdot  \left(h_r\left(t_2\right)-h_r(t_1)\right) \\
& +\sum_{\tilde{s}=m+1}^n \frac{\partial^2 G_r}{\partial u_{\tilde{s}} \partial u_s}
  \left(\tau h_r\left(t_2\right)+(1-\tau) h_r\left(t_1\right),
  \tau u_{m+1}^{(l-1)}\left(t_2, L\right)+(1-\tau) u_{m+1}^{(l-1)}(t_1, L), \cdots,\right. \\
&  \q \tau u_n^{(l-1)}\left(t_2, L\right)+(1-\tau) u_n^{(l-1)}\left(t_1, L\right) \Big)
   \cdot\left(u_{\tilde{s}}^{(l-1)}\left(t_2, L\right)-u_{\tilde{s}}^{(l-1)}\left(t_1, L\right)\right)\bigg)
   d \tau \\
&  \q \cdot z_s^{(l-1)}\left(t_1, L\right).
\end{aligned}
$$
Then we have
$$\omega\left(\delta \mid z_s^{(l)}(\cdot, 0)\right) \leq \omega\left(\delta \mid h_s^{\prime}\right)
  +C\varepsilon^2 \delta
  +\left(1+C \varepsilon \right) \theta \cdot \omega\left(\delta \mid z_r^{(l-1)}(\cdot, 0)\right)$$
and
$$\omega\left(\delta \mid z_r^{(l)}(\cdot, L)\right) \leq \omega\left(\delta \mid h_r^{\prime}\right)
  +C \varepsilon^2 \delta
  +\left(1+C \varepsilon \right) \theta \cdot \omega\left(\delta \mid z_s^{(l-1)}(\cdot, L)\right).$$
Take
\begin{equation}\label{omega}
\Omega(\delta)=\sum_{i=1}^n \frac{200}{1-\alpha} \omega\left(\delta \mid h_i^{\prime}\right)+M_4 \varepsilon \delta,
\end{equation}
where $M_4$ is a big constant which will be determined in $(\ref{m4})$,
then we get
$$\omega\left(\delta \mid z_s^{(l)}(\cdot, 0)\right) \leqslant \frac{\alpha}{8} \Omega(\delta),\q \forall s=m+1,...,n,$$
and
$$\omega\left(\delta \mid z_r^{(l)}(\cdot, L)\right) \leqslant \frac{\alpha}{8} \Omega(\delta),\q \forall r=1,...,m.$$
According to the definition $(\ref{omega})$, $\Omega\left(\delta\right)$ is continuous and
$$\lim_{\delta\rightarrow 0^+}\Omega\left(\delta\right)=0.$$
Moreover, integrate ($\ref{2.42}$) along the characteristic curve $t=t_s^{(l)}(x;t_2,x_0)$ and $t=t_s^{(l)}(x;t_1,x_0)$
from $x=0$ to $x=x_0$, we have
\begin{equation}\label{2.63}
\begin{aligned}
 & \q W_s \left(x_0\right) \left(z_s^{(l)}\left(t_2, x_0\right) - z_s^{(l)}\left(t_1,  x_0\right)\right)\\
=& z_s^{(l)}\left( t_s^{(l)}(0; t_2, x_0 ), 0\right) - z_s^{(l)}\left(t_s^{(l)}(0; t_1, x_0),  0\right)\\
 & +\left. W_s (x_0) \sum_{j\neq s} B_{s j}\left( u^{(l-1)}\right) z_j^{(l-1)}\right|_{(t_1, x_0)}^{(t_2, x_0)}
   -\left. \sum_{j\neq s} B_{s j}\left( u^{(l-1)}\right) z_j^{(l-1)}\right|_{\left( t_s^{(l)}(0; t_1, x_0 ), 0\right)}^{\left( t_s^{(l)}(0; t_2, x_0 ), 0\right) }\\
 & +\left. \int_0^{x_0}\left(-\sum_{j=1}^n \frac{\partial \mu_s}{\partial u_j}\left(u^{(l-1)}\right) z_j^{(l-1)} W_s z_s^{(l)}\right)\right|_{\left(t_s^{(l)}\left(x; t_1, x_0\right), x\right)} ^{\left(t_s^{(l)}\left(x; t_2, x_0\right), x\right)} d x \\
& +\left.\int_0^{x_0}(I_3)\right|_{\left( t_s^{(l)}(x; t_1, x_0 ), 0\right)}^{\left( t_s^{(l)}(x; t_2, x_0 ), 0\right)} d x
  +\left.\int_0^{x_0}\left(W_s \sum_{j \neq s} \tilde{g}_{s j} z_j^{(l-1)}\right)\right|
  _{\left( t_s^{(l)}(x; t_1, x_0 ), 0\right)} ^{\left( t_s^{(l)}(x; t_2, x_0 ), 0\right) } d x, \\
\end{aligned}
\end{equation}
where
\begin{align*}
 I_3= & W_s \left( \sum_{j , k=1}^n \frac{\partial B_{s j}}{\partial u_k}\left(u^{(l-1)}\right)\left(z_k^{(l-1)} \omega_j^{(l-1)} -z_j^{(l-1)} \omega_k^{(l-1)}\right)\right.\\
  & +\sum_{j, k=1}^n B_{s j}\left(u^{(l-1)}\right) \frac{\partial \mu_s}{\partial u_k}\left(u^{(l-1)}\right) z_k^{(l-1)} z_j^{(l-1)}
    +\sum_{j=1}^n \tilde{g}_{s s} B_{s j}\left(u^{(l-1)}\right) z_j^{(l-1)} \\
  &  +\sum_{j=1}^n \frac{\partial \tilde{g}_s^{N L}}{\partial u_j}\left(u^{(l-1)}\right) z_j^{(l-1)}
   \left. +K \mu_s(0) z_s^{(l-1)}\right),
    \q\q \forall s=m+1,...,n.
\end{align*}
Utilize the definition of $t_i^{(l)}(x; t_0, x_0) (i=1,...,n)$, i.e. ($\ref{2.26}$), we have
$$
\begin{aligned}
  & t_i^{(l)}\left(x; t_2, x_0\right)-t_i^{(l)}\left(x; t_1, x_0\right)\\
= & t_2-t_1+\int_{x_0}^x \left( \mu_i\left(u^{(l-1)}\left(t_i^{(l)}\left(y; t_2, x_0\right), y\right)\right)
    -\mu_i\left(u^{(l-1)}\left(t_i^{(l)}\left(y; t_1, x_0\right), y\right)\right)\right) d y \\
= & t_2-t_1+\int_{x_0}^x \int_0^1 \sum_{j=1}^n \frac{\partial \mu_i}{\partial u_j} \frac{\partial u_j^{(l-1)}}{\partial t}\left(\tau t_i^{(l)}\left(y; t_2, x_0\right)+(1-\tau) t_i^{(l)}\left(y; t_1, x_0\right), y\right) d \tau \\
& \cdot\left(t_i^{(l)}\left(y; t_2, x_0\right)-t_i^{(l)}\left(y; t_1, x_0\right)\right) d y,
\end{aligned}
$$
then applying the Gronwall's inequality, we get
\begin{equation}\label{2.64}
\left|t_i^{(l)}\left(x ; t_2, x_0\right)-t_i^{(l)}\left(x ; t_1, x_0\right)\right|
\leq (1+C \varepsilon)\left|t_2-t_1\right|
\leq (1+C \varepsilon)\delta.
\end{equation}
Due to the concavity, we have
$$
\frac{1}{1+C \varepsilon} \Omega\left(\left(1+C \varepsilon\right) \delta\right)+\frac{C \varepsilon}{1+C \varepsilon} \Omega(0)
\leq \Omega(\delta),
$$
i.e.
\begin{equation}\label{2.62}
\Omega\left((1+C \varepsilon) \delta\right) \leq(1+C \varepsilon) \Omega(\delta).
\end{equation}
Thus from ($\ref{2.63}$)-($\ref{2.62}$), and the \emph{bootstrap arguments}, we get
\begin{equation}\label{2.66}
\begin{aligned}
& \omega\left(\delta \mid z_s^{(l)}\left(\cdot, x_0\right)\right) \\
\leq    & \frac{\left(1+C \varepsilon \right) \frac{\alpha}{8} \Omega(\delta)}{W_s\left(x_0\right)}
           +\frac{C \varepsilon \Omega(\delta)+ C \varepsilon^2 \delta}{W_s\left(x_0\right)}
           +\frac{(\sum\limits_{j \neq s}\left|\tilde{g}_{s j}\right|)\left(W_s\left(x_0\right)-1\right)}{\left|\tilde{g}_{s s}\right| W_s\left(x_0\right)} \cdot
           \frac{1+C \varepsilon }{8} \Omega(\delta) \\
< & \frac{1}{8} \Omega(\delta),
\q\q  \forall s=m+1, \ldots, n .
\end{aligned}
\end{equation}
Similarly, we have
\begin{equation}\label{2.67}
\omega\left(\delta \mid z_r^{(l)}\left(\cdot, x_0\right)\right)
< \frac{1}{8} \Omega(\delta), \quad  \forall r=1, \ldots, m .
\end{equation}
Now take
\begin{equation}\label{m4}
M_4 > 1+100 \max _i\left|\tilde{g}_{i i}\right| \cdot M_1,
\end{equation}
then
$$
\left|\tilde{g}_{ii}\left(u_i^{(l)}\left(t_2,x_0\right) -u_i^{(l)}\left(t_1,x_0\right)\right)\right|
\leq \left|\tilde{g}_{ii}\right|\cdot M_1 \varepsilon \cdot \delta
\leq \frac{1}{100} M_4 \varepsilon \delta.
$$
Hence from ($\ref{2.5}$), we have
\begin{equation}\label{2.68}
\omega\left(\delta \mid \partial _x  u_i^{(l)}\left(\cdot, x_0\right)\right)
< \frac{1}{4} \Omega(\delta), \quad\q \forall i=1, \ldots, n.
\end{equation}
Thus we complete the estimation $(\ref{2.15})$.

Finally, we prove the estimation of the modulus of continuity ($\ref{2.16}$) with the help of $(\ref{2.15})$.
For $t_0 \in \mathbb{R}$, $0\leq x_1 < x_2 \leq L$ and $\left|x_2 -x_1\right|< \delta $, note $\mu_{max}\leq1$, we have
$$\left| t_i^{(l)}\left(x_2; t_0, x_1\right)-t_0 \right|\leq \delta, \q \forall i=1,...,n.$$
See Figure 1.
\begin{figure}[ht]
\begin{center}
\includegraphics[scale=0.6]{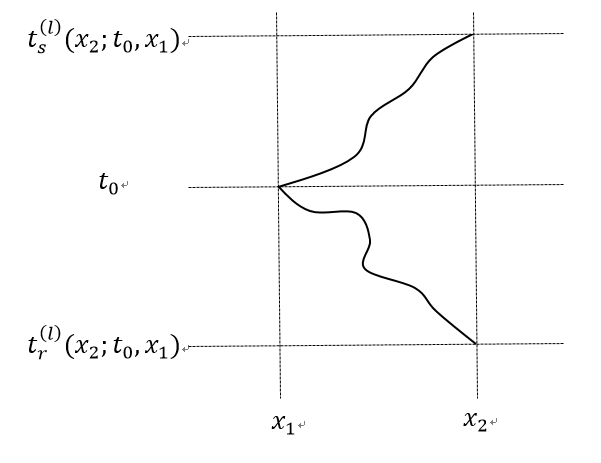}

\small{Figure 1}
\end{center}
\label{fig1}
\end{figure}

Integrate ($\ref{2.42}$) along the characteristic curve $t=t_s^{(l)}(x;t_0,x_1)$
from $x=x_1$ to $x=x_2$ and get
$$\begin{aligned}
  & W_s\left(x_2\right) z_s^{(l)}\left(t_s^{(l)}\left(x_2 ; t_0, x_1\right), x_2\right)
       -W_s\left(x_1\right) z_s^{(l)}\left(t_0, x_1\right)\\
= &\left.\sum_{j=1}^{n}\left(W_s B_{s j}\left(u^{(l-1)}\right) z_j^{(l-1)}\right)\right|
     ^{\left(t_s^{(l)}\left(x_2 ; t_0, x_1\right), x_2\right)}
     _{\left(t_0, x_1\right)}
 + \int_{x_1}^{x_2}\left(-\sum_{j=1}^n \frac{\partial \mu_s}{\partial u_j}\left(u^{(l-1)}\right) z_j^{(l-1)} W_s z_s^{(l)}\right) d x \\
  & +\int_{x_1}^{x_2} (I_3) d x+
    \left.\sum_{j \neq s} \int_{x_1}^{x_2}\left(W_s(x) \tilde{g}_{s j} z_j^{(l-1)}\right)\right|
    _{t=t_s^{(l)}\left(x ; t_0, x_1\right)} d x ,\\
\end{aligned}$$
so that
\begin{equation}\label{2.69}
\begin{aligned}
          &\left|z_s^{(l)}\left(t_s^{(l)}\left(x_2 ; t_0, x_1\right), x_2\right)-z_s^{(l)}\left(t_0, x_1\right)\right| \\
\leq      & \frac{\left|W_s\left(x_1\right)-W_s\left(x_2\right)\right|}{W_s\left(x_2\right)} \cdot M_1 \varepsilon
            +C \varepsilon^2 \delta
            +C K \varepsilon \delta+C \varepsilon \Omega(\delta) \\
          & +\frac{(\sum\limits_{j \neq s}\left|\tilde{g}_{s j}\right|)\left|W_s\left(x_1\right)-W_s\left(x_2\right)\right|}
           {\left|\tilde{g}_{s s}\right| \cdot W_s\left(x_2\right)}\cdot M_1 \varepsilon. \\
\end{aligned}
\end{equation}
Notice that
$$
\frac{\left|W_s\left(x_1\right)-W_s\left(x_2\right)\right|} { W_s\left(x_2\right)}
\leq 2\left|\tilde{g}_{ss}\right|\delta,
$$
then we have
$$ \left|z_s^{(l)}\left(t_s^{(l)}\left(x_2 ; t_0, x_1\right), x_2\right)-z_s^{(l)}\left(t_0, x_1\right)\right| \leq C \varepsilon \cdot\Omega(\delta)+\frac{1}{20} M_4 \varepsilon \delta ,$$
thus
$$\begin{aligned}
     & \left|z_s^{(l)}\left(t_0, x_1\right)-z_s^{(l)}\left(t_0, x_2\right)\right| \\
\leq & \left|z_s^{(l)}\left(t_0, x_1\right)-z_s^{(l)}\left(t_s^{(l)}\left(x_2 ; t_0, x_1\right), x_2\right)\right|+\left|z_s^{(l)}\left(t _s^{(l)}\left(x_2 ; t_0, x_1\right), x_2\right)-z_s^{(l)}\left(t_0, x_2\right)\right| \\
\leq & C \varepsilon \Omega(\delta)+\frac{1}{20} M_4 \varepsilon \delta+\frac{1}{8} \Omega(\delta) \\
\leq & \frac{1}{5} \Omega(\delta).
\end{aligned}$$
Therefore, when $\left|t_1 -t_2\right|< \delta $ and $\left|x_1 -x_2\right|< \delta $,
\begin{equation}\label{2.70}
\begin{aligned}
          & \left|z_s^{(l)}\left(t_1, x_1\right)-z_s^{(l)}\left(t_2, x_2\right)\right| \\
\leq      & \left|z_s^{(l)}\left(t_1, x_1\right)-z_s^{(l)}\left(t_1, x_2\right)\right|+\left|z_s^{(l)}\left(t_1, x_2\right)-z_s^{(l)}\left(t_2, x_2\right)\right| \\
\leq      & \frac{1}{5} \Omega(\delta)+\frac{1}{8} \Omega(\delta) \\
\leq      & \frac{13}{40} \Omega(\delta),
\q\q\q   \forall s=m+1,...,n.
\end{aligned}\end{equation}
Similarly, we have
\begin{equation}\label{2.71}
\left|z_r^{(l)}\left(t_1, x_1\right)-z_r^{(l)}\left(t_2, x_2\right)\right| \leq \frac{13}{40} \Omega(\delta),
\q\q \forall r=1,...,m.
\end{equation}
Substitute ($\ref{2.70}$) and ($\ref{2.71}$) into ($\ref{2.5}$), then we finish the prove of ($\ref{2.16}$), so as to Proposition 2.1.
\section{Stability of the time-periodic solution}
\setcounter{equation}{0}

In this section, we will prove the stability of the time-periodic solution $u^{(P)}$ obtained in Section 2.
Firstly, we give Lemma 3.1, which demonstrates the existence and uniqueness of $C^1$ solution for initial-boundary problem $(\ref{1.1})$ and $(\ref{1.16'})$-$(\ref{1.17})$ with dissipation boundary condition and $K$-weakly diagonally dominate source terms.
\\
\\
\textbf{Lemma 3.1}\emph{
Suppose assumptions $(\ref{1.15})$ and $(\ref{1.23})$ fulfil with small enough positive constant $K$, there exists a small constant $\tilde{\varepsilon}_{2}$, such that for any $\varepsilon \in (0, \tilde{\varepsilon}_{2})$, there exists $\sigma >0$, if
\begin{equation}\label{3.1}
\|u_{i0}\|_{C^{1}[0,L]}+\|h_i(t)\|_{C^{1}(\mathbb{R}^+)}\leq \sigma
\end{equation}
for $i=1,2,...,n$, the initial-boundary problem $(\ref{1.1})$ and $(\ref{1.16'})$-$(\ref{1.17})$ have an unique $C^1$ solution satisfying
\begin{equation}\label{3.1'}
\|u_{i}(t,x)\|_{C^{1}(\mathbb{R}^+ \times [0,L])}\leq C_E \varepsilon.
\end{equation}
}
\\
With the same method in \cite{Li}, we can give the proof of Lemma 3.1. Here we omit the details.

From Lemma 3.1 and Theorem 1, it is easy to see that
$$
\max_{i=1,...,n}\|u_i(t, \cdot)-u_i^{(P)}(t, \cdot)\|_{C^0}\leq 2 C_E \varepsilon,
\q \forall t\in \mathbb{R}^{+}.
$$
Then we will prove Theorem 2 inductively. Suppose for some $t_0 \geq 0$ and $N \in \mathbb{N}$, we have
\begin{equation}\label{3.2}
\max _{i=1, \ldots, n}\left\|u_i(t, \cdot)-u_i^{(P)}(t, \cdot)\right\|_{C^0} \leq C_S \varepsilon \beta^N, \quad \forall t \in\left[t_0, t_0+T_0\right]
\end{equation}
where $N=[t/ T_0]$, $C_S=2C_E$, $\beta \in (0,1)$, then we will show that
\begin{equation}\label{3.3}
\max _{i=1, \ldots, n}\left\|u_i(t, \cdot)-u_i^{(P)}(t, \cdot)\right\|_{C^0} \leq C_S \varepsilon \beta^{N+1}, \quad \forall t \in\left[t_0+T_0, t_0+2 T_0\right] .
\end{equation}
For simplicity, we denote the continuous fuction
$$
\Phi(t)=\max _i \sup _{x \in[0, L]}\left|u_i(t, x)-u_i^{(P)}(t, x)\right| .
$$
Then to prove $(\ref{3.3})$, we only need to show
\begin{equation}\label{3.5}
\Phi(t) \leq C_S \varepsilon \beta^{N+1}, \quad \forall t \in\left[t_0+T_0, \tau\right],
\end{equation}
under the hypothesis
\begin{equation}\label{3.4}
\Phi(t) \leq C_S \varepsilon \beta^N, \quad \forall t \in\left[t_0, \tau\right], \tau \in\left[t_0+T_0, t_0+2 T_0\right].
\end{equation}
To this end, we need to consider the formula of $u^{(P)}(t,x)$.
On the boundary, similarly to ($\ref{1.16}$)-($\ref{1.17}$) we have
\begin{align}
& x=0: u_s^{(P)}=G_s\left(h_s(t), u_1^{(P)}, \ldots, u_m^{(P)}\right), \q\q\q   s=m+1, \ldots, n,  \label{b1}\\
& x=L: u_r^{(P)}=G_r\left(h_r(t), u_{m+1}^{(P)}, \ldots, u_n^{(P)}\right), \q\q \ r=1, \ldots, m.  \label{b2}
\end{align}
Therefore, we have
\begin{equation}\label{3.6}
\begin{aligned}
u_r(t, L)-u_r^{(P)}(t, L)= & \sum_{s=m+1}^n\left(u_s(t, L)-u_s^{(P)}(t, L)\right)\cdot \\
& \int_0^1 \frac{\partial G_r}{\partial u_s}\left(h_r(t), \tau u_{m+1}(t, L)+(1-\tau) u_{m+1}^{(P)}(t, L),\right.\\
& \q  \left.\ldots, \tau u_n(t, L)+(1-\tau) u_n^{(P)}(t, L)\right) \mathrm{d}\tau, \quad \forall r=1, \ldots, m
\end{aligned}
\end{equation}
and
\begin{equation}\label{3.7}
\begin{aligned}
u_s(t, 0)-u_s^{(P)}(t, 0)= & \sum_{r=1}^m\left(u_r(t, 0)-u_r^{(P)}(t, 0)\right)\cdot \\
& \int_0^1 \frac{\partial G_s}{\partial u_r}\left(h_s(t), \tau u_1(t, 0)+(1-\tau) u_1^{(P)}(t, 0),\right.\\
&\q \left.\ldots, \tau u_m(t, 0)+(1-\tau) u_m^{(P)}(t, 0)\right) \mathrm{d} \tau, \quad \forall s=m+1, \ldots, n.
\end{aligned}
\end{equation}
Consequently, from ($\ref{1.23}$) and  ($\ref{3.4}$), we have
\begin{equation}\label{3.8}
\max _{r=1, \ldots, m} \sup _{t \in\left[t_0, \tau\right]}\left|u_r(t, L)-u_r^{(P)}(t, L)\right|
\leq \theta (1+\varepsilon) C_S \varepsilon \beta^N,
\end{equation}
and
\begin{equation}\label{3.9}
\max _{s=m+1, \ldots, n} \sup _{t \in\left[t_0, \tau\right]}\left|u_s(t, 0)-u_s^{(P)}(t, 0)\right|
\leq \theta (1+\varepsilon)  C_S \varepsilon \beta^N.
\end{equation}

In the domain, similarly to ($\ref{2.8'}$), the time-periodic solution $u= u^{(P)}(t,x)$ satisfy the following form
\begin{align}
  \left( \partial_x u_i^{(P)}+\mu_i\left(u^{(P)}\right) \partial_t u_i^{(P)}\right)
   & -\tilde{g}_{i i} u_i^{(P)}
 =  \sum_{j=1}^n B_{i j}\left(u^{(P)}\right)\left(\partial_x +\mu_i\left(u^{(P)}\right) \partial_t \right) u_j^{(P)} \nonumber\\
   & +\sum_{j \neq i} \tilde{g}_{i j} u_j^{(P)}+K \mu_i(0) u_i^{(P)}+\tilde{g}_i^{N L}\left(u^{(P)}\right),
    \q i=1,...,n.\label{3.11}
\end{align}
Then multiply ($\ref{2.8'}$)($\ref{3.11}$) with $W_{r}(x)$ for $r=1,...,m$ and combine the results together, we have
\begin{equation}\label{3.12}
\begin{aligned}
  & \left(\partial_x + \mu_r(u)\partial_t \right)\left(W_r\left(u_r-u_r^{(P)}\right)\right) \\
= & \left(\partial_x + \mu_r(u)\partial_t \right)\left(W_r u_r \right)
    - \left(\partial_x + \mu_r(u^{(P)})\partial_t \right)\left(W_r u_r^{(P)} \right)\\
    & -W_r\cdot\left(\mu_r\left(u\right)-\mu_r\left(u^{(P)}\right)\right) \partial_t u_r^{(P)}\\
= & \left(\partial_x+\mu_r\left(u\right) \partial_t\right)\left(\sum_{j=1}^n W_r B_{r j}\left(u\right)\left(u_j-u_j^{(P)}\right)\right)+I_4 \\
& +W_r\left(\beta \sum_{j \neq r} \tilde{g}_{r j}\left(u_j-u_j^{(P)}\right)+(1-\beta) \sum_{j \neq r} \tilde{g}_{r j}\left(u_j-u_j^{(P)}\right)\right),
\q \forall r=1,...,m,
\end{aligned}
\end{equation}
where
\begin{equation}\label{I4}
\begin{aligned}
I_4 = & W_r\left(-\sum_{j=1}^n \partial_t u_r \int_0^1 \frac{\partial \mu_r}{\partial u_j}\left(\tau u+(1-\tau) u^{(P)}\right) d \tau \cdot\left(u_j -u_j^{(P)}\right)\right. \\
  & -\sum_{j,k=1}^n \frac{\partial B_{r j}}{\partial u_k}\left(u\right)\left(\partial_x +\mu_r\left(u\right) \partial_t \right) u_k \cdot \left(u_j-u_j^{(P)}\right)
   +\sum_{j=1}^{n} B_{r j}\left(u\right) \tilde{g}_{r r}\left(u_j-u_j^{(P)}\right) \\
  & +\sum_{j,k=1}^{n}\left(\partial_x u_j^{(P)}+\mu_r\left(u\right) \partial_t u_j^{(P)}\right) \int_0^1 \frac{\partial B_{r j}}{\partial u_k}\left(\tau u+(1-\tau) u^{(P)}\right) d \tau \cdot\left(u_k-u_k^{(P)}\right) \\
  & +\sum_{j, k=1}^n  B_{r j}\left(u\right) \partial_t u_j^{(P)} \int_0^1 \frac{\partial \mu_r}{\partial u_k}\left(\tau u+(1-\tau) u^{(P)}\right) d \tau \cdot\left(u_k-u_k^{(P)}\right) \\
  & +\sum_{j=1}^n  \int_0^1 \frac{\partial \tilde{g}_r^{N L}}{\partial u_j}\left(\tau u_j+(1-\tau) u_j^{(P)}\right) d \tau \cdot\left(u_j-u_j^{(P)}\right) \\
  &\left.  +K \mu_r(0)\left(u_r-u_r^{(P)}\right) \right),
  \q\q\q\q\q  \forall r=1,...,m.\\
\end{aligned}
\end{equation}
Similarly, we have
\begin{equation}\label{3.13}
\begin{aligned}
  & \left(\partial_x + \mu_s(u)\partial_t \right)\left(W_s\left(u_s-u_s^{(P)}\right)\right) \\
= & \left(\partial_x + \mu_s(u)\partial_t \right)\left(W_s u_s \right)
    - \left(\partial_x + \mu_s(u^{(P)})\partial_t \right)\left(W_s u_s^{(P)} \right)\\
   & -W_s\cdot\left(\mu_s\left(u\right)-\mu_s\left(u^{(P)}\right)\right) \partial_t u_s^{(P)}\\
= & \left(\partial_x+\mu_s\left(u\right) \partial_t\right)\left(\sum_{j=1}^n W_s B_{s j} \left(u\right) \left(u_j-u_j^{(P)}\right) \right) +I_5 \\
& +W_s\left(\beta \sum_{j \neq s} \tilde{g}_{s j}\left(u_j-u_j^{(P)}\right)+(1-\beta) \sum_{j \neq s} \tilde{g}_{s j}\left(u_j-u_j^{(P)}\right)\right),
 \q \forall s=m+1,...,n,
\end{aligned}
\end{equation}
with
\begin{equation}\label{I5}
\begin{aligned}
I_5 = & W_s\left(-\sum_{j=1}^n \partial_t u_s \int_0^1 \frac{\partial \mu_s}{\partial u_j}\left(\tau u+(1-\tau) u^{(P)}\right) d \tau \cdot\left(u_j -u_j^{(P)}\right)\right. \\
  & -\sum_{j,k=1}^n \frac{\partial B_{s j}}{\partial u_k}\left(u\right)\left(\partial_x +\mu_s\left(u\right) \partial_t \right) u_k \cdot \left(u_j-u_j^{(P)}\right)
   +\sum_{j=1}^{n} B_{s j}\left(u\right) \tilde{g}_{s s}\left(u_j-u_j^{(P)}\right) \\
  & +\sum_{j,k=1}^{n}\left(\partial_x u_j^{(P)}+\mu_s\left(u\right) \partial_t u_j^{(P)}\right) \int_0^1 \frac{\partial B_{s j}}{\partial u_k} \left(\tau u+(1-\tau) u^{(P)}\right) d \tau \cdot\left(u_k-u_k^{(P)}\right) \\
  & +\sum_{j, k=1}^n  B_{s j}\left(u\right) \partial_t u_j^{(P)} \int_0^1 \frac{\partial \mu_s}{\partial u_k}\left(\tau u+(1-\tau) u^{(P)}\right) d \tau \cdot\left(u_k-u_k^{(P)}\right) \\
  & +\sum_{j=1}^n  \int_0^1 \frac{\partial \tilde{g}_s^{N L}}{\partial u_j}\left(\tau u_j+(1-\tau) u_j^{(P)}\right) d \tau \cdot\left(u_j-u_j^{(P)}\right) \\
  &\left.  +K \mu_s(0)\left(u_s-u_s^{(P)}\right) \right),
  \q\q\q\q\q  \forall s=m+1,...,n.\\
\end{aligned}
\end{equation}

Integrate $(\ref{3.13})$ along the characteristic curve $t= t(x; t_0,x_0)$ from $x = 0$ to $x = x_0$ and notice ($\ref{beta}$), we have
\begin{equation}\label{3.14}
\begin{aligned}
   & \left|u_s\left(t_0, x_0\right)-u_s^{(P)}\left(t_0, x_0\right)\right|\\
 \leq & \frac{\theta (1+ C \varepsilon)}{W_s (x_0)} C_S \varepsilon\beta^{N}
        +\frac{C (\varepsilon + K)}{W_s (x_0)} C_S \varepsilon\beta^{N}
        +\frac{\left( \sum\limits_{j\neq s}\left|\tilde{g}_{sj}\right| \right)\left( W_s(x_0)-1 \right)}{\left| \tilde{g}_{ss} \right|W_s (x_0)} C_S \varepsilon\beta^{N+1} \\
      & +\frac{\left( \sum\limits_{j\neq s}\left|\tilde{g}_{sj}\right| \right)\left( W_s(x_0)-1 \right)}{\left| \tilde{g}_{ss} \right|W_s (x_0)} \left(1-\beta\right)C_S \varepsilon\beta^{N} \\
 \leq & C_S \varepsilon\beta^{N+1}+\frac{C_S \varepsilon\beta^{N}}{W_s\left(x_0\right)}
      \bigg(\theta\left(1+C \varepsilon \right)
      +C\left(\varepsilon+K\right)+\left(1-\beta\right)W_s\left(x_0\right)-1 \bigg)\\
 \leq & C_S \varepsilon \beta^{N+1}.
\end{aligned}
\end{equation}
Similarly, integrate $(\ref{3.12})$ along the characteristic curve $t= t(x; t_0,x_0)$ from $x = L$ to $x = x_0$, we have
\begin{equation}
\left|u_r\left(t_0, x_0\right)-u_r^{(P)}\left(t_0, x_0\right)\right|  \leq C_S \varepsilon \beta^{N+1},
\q\q \forall r=1,...,m.
\end{equation}
Thus
\begin{equation}\label{3.15}
\Phi(t) \leq C_S \varepsilon \beta^{N+1},
\end{equation}
and we complete the proof of Theorem 2.
\section{Regularity of the time-periodic solution}
\setcounter{equation}{0}

In this section, we will give the proof of Theorem 4, i.e., the regularity of the time-periodic solution $u^{(P)}$ with the $W^{2,\infty}$ regularity of boundary functions $h_i (t)$.
First, we use the iterative scheme $(\ref{2.5})$-$(\ref{2.12})$ and give the following proposition:
 \\\\
\textbf{Proposition 4.1.} \emph{Suppose the assumptions of Theorem 4 fulfil, there exists a large enough constant $M_R >0$, such that the solutions of the iteration system $(\ref{2.5})$-$(\ref{2.12})$ satisfy
\begin{equation}\label{4.2}
\max_{i=1,...,n}\left\{ \|\partial_t^2 u^{(l)}\|_{L^\infty}, \|\partial_t \partial_x u^{(l)}\|_{L^\infty},
\|\partial_x^2 u^{(l)}\|_{L^\infty},\right\} \leq M_R < +\infty,
\end{equation}
under the hypothesis
\begin{equation}\label{4.1}
\max_{i=1,...,n}\left\{ \|\partial_t^2 u^{(l-1)}\|_{L^\infty}, \|\partial_t \partial_x u^{(l-1)}\|_{L^\infty},
\|\partial_x^2 u^{(l-1)}\|_{L^\infty},\right\} \leq M_R < +\infty.
\end{equation}
}

Once Proposition 4.1 is proved, we get the uniform $W^{2,\infty}$ bounded of the sequence $\{u^{(l)}\}_{l=1}^{\infty}$ from $(\ref{4.2})$. Consequently the weak$\ast$ convergence arrives.
Based on the strong convergence of $\{u^{(l)}\}_{l=1}^{\infty}$ in Theorem 1, then we show the $W^{2,\infty}$ regularity of $u^{(P)}$.
\\\\
\textbf{Proof :}
Actually, since Proposition 2.1, for each $l \in \mathbb{Z}_{+}$, we have  $(\ref{2.17})$-$(\ref{2.16})$, and
\begin{equation}\label{4.3}
\left\|u^{(l)}\right\|_{C^{1}} \leq M_1 \varepsilon,
\q \left\|u^{(l-1)}\right\|_{C^{1}} \leq M_1 \varepsilon.
\end{equation}
Denote
\begin{equation}\label{4.4}
\phi_i^{(l)}= \partial_t z_i^{(l)}= \partial_t^2 u_i^{(l)}, \q i=1,...,n; \ l \in \mathbb{Z}_{+},
\end{equation}
then take the temporal derivative on $(\ref{2.36})$-$(\ref{2.37})$, we have
\begin{equation}\label{4.5}
\begin{aligned}
\quad x=0: \phi_s^{(l)}
&=\frac{\partial G_s}{\partial h_s}\left(h_s, u_1^{(l-1)},..., u_m^{(l-1)}\right) h_s^{''}(t)
   + \frac{\partial^2 G_s}{\partial h_s^2}\left(h_s, u_1^{(l-1)},..., u_m^{(l-1)}\right) (h_s^{'}(t))^2\\
& + 2 \sum_{r=1}^m \frac{\partial^2 G_s}{\partial h_s \partial u_r}\left(h_s, u_1^{(l-1)}, \ldots, u_m^{(l-1)}\right) z_r^{(l-1)}h^{'}_s (t)\\
& +\sum_{r=1}^m \frac{\partial G_s}{\partial u_r}\left(h_s, u_1^{(l-1)}, \ldots, u_m^{(l-1)}\right) \phi_r^{(l-1)}\\
&  +\sum_{r,r'=1}^m \frac{\partial^2 G_s}{\partial u_r \partial u_{r'}}\left(h_s, u_1^{(l-1)}, \ldots, u_m^{(l-1)}\right) z_r^{(l-1)} z_{r'} ^{(l-1)} ,
\q  \forall s=m+1,...,n,
\end{aligned}\end{equation}
and
\begin{equation}\label{4.6}
\begin{aligned}
\quad x=L: \phi_r^{(l)}
&=\frac{\partial G_r}{\partial h_r}\left(h_r, u_{m+1}^{(l-1)},..., u_n^{(l-1)}\right) h_r^{''}(t)
   + \frac{\partial^2 G_r}{\partial h_r^2}\left(h_r, u_{m+1}^{(l-1)},..., u_n^{(l-1)}\right) (h_r^{'}(t))^2\\
& + 2 \sum_{s=m+1}^n \frac{\partial^2 G_r}{\partial h_r \partial u_s}\left(h_r, u_{m+1}^{(l-1)}, \ldots, u_n^{(l-1)}\right) z_s^{(l-1)}h^{'}_r (t)\\
& +\sum_{s=m+1}^n \frac{\partial G_r}{\partial u_s}\left(h_r, u_{m+1}^{(l-1)}, \ldots, u_n^{(l-1)}\right) \phi_s^{(l-1)}\\
&  +\sum_{s, s'=m+1}^n \frac{\partial^2 G_r}{\partial u_s \partial u_{s'}}\left(h_r, u_{m+1}^{(l-1)}, \ldots, u_n^{(l-1)}\right) z_s^{(l-1)} z_{s'} ^{(l-1)} ,
\q  \forall r=1,...,m.
\end{aligned}\end{equation}
By assumptions $(\ref{1.20})$, $(\ref{1.23})$, $(\ref{4.1})$, and the boundary regularity $(\ref{1.29})$,
we can easily arrive
\begin{equation}\label{4.7}
\sup_t | \phi_r^{(l)}(t, L)| \leq (\frac{1}{2} + C \varepsilon)M_0 + (\frac{1}{2}+ C \varepsilon)M_R  + C \varepsilon^2
                             \leq \lambda M_R,
\q \forall r=1,...,m
\end{equation}
for small $\lambda\in (0,1)$ and lager enough constant $M_R$.
Similarly,
\begin{equation}\label{4.8}
\sup_t | \phi_s^{(l)}(t, 0)| \leq (\frac{1}{2} + C \varepsilon)M_0 + (\frac{1}{2}+ C \varepsilon)M_R  + C \varepsilon^2
                             \leq \lambda M_R,
\q \forall s=m+1,...,n.
\end{equation}

In the domain $\{(t,x)| t \in \mathbb{R} ,0<x< L\}$, we also take the temporal derivative on $(\ref{2.40})$ to get
\begin{equation}\label{4.9}
\begin{aligned}
  & \partial_x \phi_i^{(l)}+\mu_i\left(u^{(l-1)}\right) \partial_t \phi_i^{(l)}-\tilde{g}_{i i} \phi_i^{(l)}\\
= & -2\sum_{j=1} ^n  \displaystyle\frac{\partial \mu_i}{\partial u_j}\left(u^{(l-1)}\right) z_j^{(l-1)} \phi_i^{(l)}
    -\sum_{j,k =1} ^n  \displaystyle\frac{\partial^2 \mu_i}{\partial u_j \partial u_k }\left(u^{(l-1)}\right) z_k^{(l-1)} z_j^{(l-1)} z_i^{(l)}\\
    & -\sum_{j=1} ^n  \displaystyle\frac{\partial \mu_i}{\partial u_j}\left(u^{(l-1)}\right) \phi_j^{(l-1)} z_i^{(l)}\\
        &  +\sum_{j, k, p=1}^n \frac{\partial^2 B_{i j}}{\partial u_k \partial u_p }\left(u^{(l-1)}\right) z_p^{(l-1)} z_k^{(l-1)}\left(\omega_j^{(l-1)}+\mu_i(u^{(l-1)}) z_j^{(l-1)}\right)\\
    & +\sum_{j, k=1}^n \frac{\partial B_{i j}}{\partial u_k}\left(u^{(l-1)}\right) \phi_k^{(l-1)}\left(\omega_j^{(l-1)}+\mu_i\left(u^{(l-1)}\right) z_j^{(l-1)}\right)\\
    & +\sum_{j, k=1}^n \frac{\partial B_{i j}}{\partial u_k}\left(u^{(l-1)}\right) z_k^{(l-1)}\left(\partial_x z_j^{(l-1)}+\mu_i\left(u^{(l-1)}\right) \partial_t z_j^{(l-1)}\right)\\
    & +\sum_{j, k, p=1}^n \frac{\partial B_{i j}}{\partial u_k}\left(u^{(l-1)}\right) z_k^{(l-1)} \displaystyle\frac{\partial \mu_i}{\partial u_p} \left(u^{(l-1)}\right)  z_p^{(l-1)} z_j^{(l-1)}\\
    & +\sum_{j, k, p=1}^n \frac{\partial B_{i j}}{\partial u_p}\left(u^{(l-1)}\right) \displaystyle\frac{\partial \mu_i}{\partial u_p}\left(u^{(l-1)}\right)  z_j^{(l-1)} z_p^{(l-1)}\\
    & +\sum_{j, k, p=1} ^n  B_{i j}\left(u^{(l-1)}\right) \frac{\partial^2 \mu_i}{\partial u_k \partial u_p}\left(u^{(l-1)}\right) z_k^{(l-1)} z_j^{(l-1)} z_p^{(l-1)}\\
    & +\sum_{j, k=1} ^n  B_{i j}\left(u^{(l-1)}\right) \frac{\partial \mu_i}{\partial u_k }\left(u^{(l-1)}\right) \left(\phi_k^{(l-1)} z_j^{(l-1)} + z_k^{(l-1)} \phi_j^{(l-1)} \right)\\
\end{aligned}
\end{equation}
$$
\begin{aligned}
    & +\sum_{j=1}^n \left(\partial_x+\mu_i\left(u^{(l-1)}\right) \partial_t\right)\left(B_{i j}\left(u^{(l-1)}\right) \phi_j^{(l-1)}
    + \displaystyle\frac{ \partial B_{i j}}{\partial u_k} \left(u^{(l-1)}\right) z_j^{(l-1)} z_k^{(l-1)}\right)\\
    & +\sum_{j, k=1}^n \left(\partial_x+\displaystyle\frac{\partial \mu_i}{\partial u_k}\left(u^{(l-1)}\right) \partial_t\right)\left(B_{i j}\left(u^{(l-1)}\right) z_j^{(l-1)}\right)\\
    & -\sum_{j, k, p=1}^n \frac{\partial^2 B_{i j}}{\partial u_k \partial u_p }\left(u^{(l-1)}\right) z_p^{(l-1)} z_j^{(l-1)}\left(\omega_k^{(l-1)}+\mu_i(u^{(l-1)}) z_k^{(l-1)}\right)\\
      & -\sum_{j, k=1}^n \frac{\partial B_{i j}}{\partial u_k}\left(u^{(l-1)}\right) \phi_j^{(l-1)}\left(\omega_k^{(l-1)}+\mu_i\left(u^{(l-1)}\right) z_k^{(l-1)}\right)\\
     & -\sum_{j, k=1}^n \frac{\partial B_{i j}}{\partial u_k}\left(u^{(l-1)}\right) z_j^{(l-1)}\left(\partial_x z_k^{(l-1)}+\mu_i\left(u^{(l-1)}\right) \partial_t z_k^{(l-1)}\right)\\
     & -\sum_{j, k, p=1}^n \frac{\partial B_{i j}}{\partial u_k}\left(u^{(l-1)}\right) z_j^{(l-1)} \displaystyle\frac{\partial \mu_i}{\partial u_p}\left(u^{(l-1)}\right)  z_k^{(l-1)} z_p^{(l-1)}\\
    & +\sum_{j, k=1}^n \frac{\partial^2 \tilde{g}_i^{N L}}{\partial u_j \partial u_k}\left(u^{(l-1)}\right) z_k^{(l-1)} z_j^{(l-1)}
      +\sum_{j=1}^n \frac{\partial \tilde{g}_i^{N L}}{\partial u_j}\left(u^{(l-1)}\right) \phi_j^{(l-1)}\\
    & +K \mu_i(0) \phi_i^{(l-1)}+\sum_{j \neq i} \tilde{g}_{i j} \phi_j^{(l-1)},
\q\q  \forall i=1,...,n.
\end{aligned}
$$
For $r=1,...,m$, we rewrite $(\ref{4.9})$ with $W_r(x)$ and integrate along the characteristic curve $t= t_r^{(l)}(x, t_0, x_0)$
from $x=L$ to $x=x_0$ to get
\begin{equation}\label{4.10}
\begin{aligned}
 \| \phi_r^{(l)} \|_{L^{\infty}}
 \leq & \frac{\lambda M_R}{W_r (x_0)}
        +\frac{C (\varepsilon + K)}{W_r (x_0)} \cdot M_R
        +\frac{\left( \sum\limits_{j\neq r}\left|\tilde{g}_{rj}\right| \right)\left( W_r(x_0)-1 \right)}{\left| \tilde{g}_{rr} \right|W_r (x_0)}\cdot M_R \\
 \leq & M_R+\frac{M_R}{W_r\left(x_0\right)}
      \bigg(\lambda+C\left(\varepsilon+K\right)-1 \bigg)\\
 \leq & M_R,
\end{aligned}
\end{equation}
where we have used $(\ref{4.1})$, $(\ref{4.3})$, $(\ref{4.7})$, and the smallness of $\v$ and $K$.
Then for $s=m+1,...,n$, we can similarly have
\begin{equation}\label{4.11}
\begin{aligned}
 \| \phi_s^{(l)} \|_{L^{\infty}}
 \leq & \frac{\lambda M_R}{W_s (x_0)}
        +\frac{C (\varepsilon + K)}{W_s (x_0)} \cdot M_R
        +\frac{\left( \sum\limits_{j\neq s}\left|\tilde{g}_{sj}\right| \right)\left( W_s(x_0)-1 \right)}{\left| \tilde{g}_{ss} \right|W_s (x_0)}\cdot M_R \\
 \leq & M_R+\frac{M_R}{W_s\left(x_0\right)}
      \bigg(\lambda+C\left(\varepsilon+K\right)-1 \bigg)\\
 \leq & M_R.
\end{aligned}
\end{equation}

Next, utilize $(\ref{2.41})$, the assumptions $(\ref{1.13})$,  $(\ref{2.3})$, $(\ref{4.10})$ and $(\ref{4.11})$ lead us to get
\begin{equation}\label{4.12}
\| \partial_x \partial_t u_i^{(l)} \|_{L^{\infty}} \leq  M_R, \q \forall i=1,...,n.
\end{equation}
Moreover, take  the spatial derivative on $(\ref{2.5})$, we get
\begin{equation}\label{4.13}
\begin{aligned}
   &\partial_x^2 u_i^{(l)}+\mu_i\left(u^{(l-1)}\right)\partial_x \partial_t u_i^{(l)}-\tilde{g}_{i i} \partial_xu_i^{(l)}\\
 = & -\sum_{j=1} ^n  \displaystyle\frac{\partial \mu_i}{\partial u_j}\partial_t u_i^{(l)}
     +\sum_{j, k=1}^n \displaystyle\frac{\partial B_{i j}}{\partial u_k} \left(u^{(l-1)}\right) \left(\partial_x u_j^{(l-1)}+\mu_i\left(u^{(l-1)}\right) \partial_t u_j^{(l-1)}\right) \partial_x u_k^{(l-1)} \\
   & +\sum_{j=1}^n B_{i j}\left(u^{(l-1)}\right)\left(\partial_x^2 u_j^{(l-1)}+\mu_i\left(u^{(l-1)}\right) \partial_x \partial_t u_j^{(l-1)}\right) \\
   & +\sum_{j, k=1}^n B_{i j}\left(u^{(l-1)}\right)\left(\partial_x u_j^{(l-1)}+\displaystyle\frac{\partial \mu_i}{\partial u_k } \left(u^{(l-1)}\right) \partial_x u_k^{(l-1)} \partial_t u_j^{(l-1)}\right) \\
   & +\sum_{j\neq i} \tilde{g}_{i j} \partial_x u_j^{(l-1)}+K \mu_i(0) \partial_x u_i^{(l-1)}+\sum_{j=1}^{n} \displaystyle\frac{\partial \tilde{g}_{i}^{N L}}{\partial u_j}\left(u^{(l-1)}\right)\partial_x u_j^{(l-1)},
 \q\q \forall i=1,...,n.
\end{aligned}
\end{equation}
Combine $(\ref{1.13})$, $(\ref{2.3})$, $(\ref{4.1})$, $(\ref{4.3})$, and $(\ref{4.12})$, we have
\begin{equation}\label{4.14}
\| \partial_x^2  u_i^{(l)} \|_{L^{\infty}} \leq  M_R, \q \forall i=1,...,n.
\end{equation}
Hence we finish the proof of Proposition 4.1.
\section{Boundary stabilization around the time-periodic solution}
\setcounter{equation}{0}

In this section, we will prove Theorem 5.
Actually, same as in \cite{Li,Qin,Zhao}, we can view the boundary condition $(\ref{1.16})$-$(\ref{1.17})$ as a feedback boundary control while all the $h_i(t)\equiv 0$ for $i=1,...,n$.
When the dissipative hypothesis $(\ref{1.23})$ fulfil, the weak dissipation system $(\ref{1.1})$ can be stabilized by the feedback control near the constant equilibrium $u=0$.
In this paper, we will show the stabilization property for the system around the corresponding time-periodic solution $u^{(P)}$.

In fact, we have already obtain the $C^0$ convergence results in Theorem 2:
\begin{equation}\label{5.1}
\max _{i=1, \ldots, n}\left\|u_i(t, \cdot)-u_i^{(P)}(t, \cdot)\right\|_{C^0} \leq C_S \varepsilon \beta^N,
\quad \forall t \in\left[N T_0, (N+1)T_0\right),
\end{equation}
and further we have
\begin{equation}\label{5.2}
\max _{i=1, \ldots, n}\left\|u_i(t, \cdot)-u_i^{(P)}(t, \cdot)\right\|_{C^0} \leq C_S \varepsilon \beta^{N+1},
\quad \forall t \in\left[(N+1) T_0, (N+2)T_0\right)
\end{equation}
for each $N \in \mathbb{N}_+$.
What's more, according to Theorem 1, Lemma 3.1 and Theorem 4, we get
\begin{equation}\label{5.3}
\|u\|_{C^1} \leq C_E \varepsilon,
\quad\left\|u^{(P)}\right\|_{C^1} \leq C_E \varepsilon,
\quad\left\|u^{(P)}\right\|_{W^{2, \infty}} \leq M_R.
\end{equation}
Since continuity, we will prove the $C^0$ estimates of the first derivatives by an iteration, that is, suppose
\begin{equation}
\begin{aligned}\label{5.4}
      & \max _{i=1, \ldots, n}\left\{\left\|\partial_t u_i(t, \cdot)-\partial_t u_i^{(P)}(t, \cdot)\right\|_{C^0},\left\|\partial_x u_i(t, \cdot)-\partial_x u_i^{(P)}(t, \cdot)\right\|_{C^0}\right\}\\
 \leq & \widetilde{C}_S \beta^{N} \varepsilon,
 \q\q  \forall t \in\left[N T_0, \tau\right],
\end{aligned}
\end{equation}
for each $N \in \mathbb{N}_+$ and $\tau \in \left[(N+1) T_0, (N+2)T_0\right)$,
then we will show that
\begin{equation}
\begin{aligned}\label{5.5}
      & \max _{i=1, \ldots, n}\left\{\left\|\partial_t u_i(t, \cdot)-\partial_t u_i^{(P)}(t, \cdot)\right\|_{C^0},\left\|\partial_x u_i(t, \cdot)-\partial_x u_i^{(P)}(t, \cdot)\right\|_{C^0}\right\}\\
 \leq &  \widetilde{C}_S  \beta^{N+1} \varepsilon,
 \q\q  \forall t \in\left[(N+1) T_0, \tau\right].
\end{aligned}
\end{equation}
Similar steps to other sections, we first consider boundary estimates.
Take the temporal derivative on boundary conditions $(\ref{1.16})$-$(\ref{1.17})$, we have
\begin{align}
x=0: \partial_t u_s =& z_s=  h'_s(t) \frac{\partial G_s}{\partial h_s}\left(h_s(t), u_1, \ldots, u_m\right)  \nonumber\\
& +\sum_{r=1}^m z_r \frac{\partial G_s}{\partial u_r}\left(h_s(t), u_1, \ldots, u_m\right), \quad s=m+1, \ldots, n, \label{5.7}\\
x=L: \partial_t u_r =& z_r=  h'_r(t) \frac{\partial G_r}{\partial h_r}\left(h_r(t), u_{m+1}, \ldots, u_n\right)  \nonumber\\
& +\sum_{s=m+1}^n z_s \frac{\partial G_r}{\partial u_s}\left(h_r(t), u_{m+1}, \ldots, u_n\right), \quad r=1, \ldots, m\label{5.8}
\end{align}
and
\begin{align}
x=0: \partial_t u_s^{(P)} =& z_s^{(P)}=  h'_s(t) \frac{\partial G_s}{\partial h_s}\left(h_s(t), u_1^{(P)}, \ldots, u_m^{(P)}\right)  \nonumber\\
& +\sum_{r=1}^m z_r^{(P)} \frac{\partial G_s}{\partial u_r}\left(h_s(t), u_1^{(P)}, \ldots, u_m^{(P)}\right), \quad s=m+1, \ldots, n,\label{5.9}  \\
x=L: \partial_t u_r^{(P)} =& z_r^{(P)}=  h'_r(t) \frac{\partial G_r}{\partial h_r}\left(h_r(t), u_{m+1}^{(P)}, \ldots, u_n^{(P)}\right)  \nonumber\\
& +\sum_{s=m+1}^n z_s^{(P)} \frac{\partial G_r}{\partial u_s}\left(h_r(t), u_{m+1}^{(P)}, \ldots, u_n^{(P)}\right), \quad r=1, \ldots, m \label{5.10} .
\end{align}
For $s=m+1,...,n$, we have

\begin{equation}\label{5.11}
\begin{aligned}
& z_s (t, 0)  -  z_s^{(P)}(t, 0)
  =  h'_s(t) \sum_{r=1}^m\left(u_r(t, 0)-u_r^{(P)}(t, 0)\right) \cdot \\
& \int_0^1 \frac{\partial^2 G_s}{\partial u_r \partial h_s}\left(h_s(t), \gamma u_1(t, 0)+(1-\gamma) u_1^{(P)}(t, 0), \ldots, \gamma u_m(t, 0)+(1-\gamma) u_m^{(P)}(t, 0)\right) \mathrm{d} \gamma \quad \\
& +\sum_{r=1}^m\left(z_r(t, 0)-z_r^{(P)}(t, 0)\right) \frac{\partial G_s}{\partial u_r}\left(h_s(t), u_1(t, 0), \ldots, u_m(t, 0)\right) \\
& +\sum_{r, r^{\prime}=1}^m z_r^{(P)}(t, 0)\left(u_{r^{\prime}}(t, 0)-u_{r^{\prime}}^{(P)}(t, 0)\right) . \\
& \int_0^1 \frac{\partial^2 G_s}{\partial u_{r^{\prime}} \partial u_r}\left(h_s(t), \gamma u_1(t, 0)+(1-\gamma) u_1^{(P)}(t, 0), \ldots, \gamma u_m(t, 0)+(1-\gamma) u_m^{(P)}(t, 0)\right) \mathrm{d} \gamma.
\end{aligned}
\end{equation}
By $(\ref{5.1})$-$(\ref{5.4})$, on the boundary $x=0$ we have
\begin{equation}\label{5.12}
\begin{aligned}
\sup _{t \in\left[N T_0, \tau\right]}\left|z_s(t, 0)-z_s^{(P)}(t, 0)\right|
& \leq  \left(\frac{1}{2}+ C \varepsilon\right) \widetilde{C}_S \beta^{N} \varepsilon
        + C \varepsilon C_S \beta^{N+1} \varepsilon + C\varepsilon^2\\
& \leq \beta^{*} \widetilde{C}_S \beta^{N} \varepsilon\leq  \widetilde{C}_S \beta^{N+1} \varepsilon,
\end{aligned}
\end{equation}
where we have used the facts $1-\beta \ll 1$ and  $\frac{1}{2}+C \varepsilon < \beta^{*}< \beta $.
Similarly, for $r=1, \ldots, m$, we have
\begin{equation}\label{5.13}
\sup _{t \in\left[N T_0, \tau\right]}\left|z_r(t, L)-z_r^{(P)}(t, L)\right|
\leq \beta^{*} \widetilde{C}_S \beta^{N} \varepsilon\leq  \widetilde{C}_S \beta^{N+1} \varepsilon.
\end{equation}

In the domain $\{(t,x)|t \in \mathbb{R} ,0<x< L\}$, we consider the formulas of $z_i$ and $z_i^{(P)}$ for $i=1,...,n$.
Take the temporal derivative to $(\ref{2.8'})$  and $(\ref{3.11})$ to get
\begin{equation}\label{5.15}
\begin{aligned}
  & \partial_x z_i+\mu_i\left(u\right) \partial_t z_i-\tilde{g}_{i i} z_i\\
= & -\sum_{j=1} ^n  \displaystyle\frac{\partial \mu_i}{\partial u_j}\left(u\right) z_j z_i
    +\sum_{j=1} ^n  B_{i j}\left(u\right) \left( \partial_x z_j + \mu_i\left(u\right)\partial_t z_j \right)\\
  &  +\sum_{j, k=1} ^n  B_{i j}\left(u\right) \frac{\partial \mu_i}{\partial u_k}\left(u\right) z_k z_j
     +\sum_{j, k=1}^n \frac{\partial B_{i j}}{\partial u_k}\left(u\right) z_k\left(\omega_j+\mu_i\left(u\right) z_j\right)\\
  & +\sum_{j \neq i} \tilde{g}_{i j} z_j +K \mu_i(0) z_i
    +\sum_{j=1}^n \frac{\partial \tilde{g}_i^{N L}}{\partial u_j}\left(u\right) z_j,
\q\q  \forall i=1,...,n
\end{aligned}\end{equation}
and
$$
\begin{aligned}
  & \partial_x z_i^{(P)}+\mu_i\left(u^{(P)}\right) \partial_t z_i^{(P)}-\tilde{g}_{i i} z_i^{(P)}\\
= & -\sum_{j=1} ^n  \displaystyle\frac{\partial \mu_i}{\partial u_j}\left(u^{(P)}\right) z_j^{(P)} z_i^{(P)}
      +\sum_{j=1} ^n  B_{i j}\left(u^{(P)}\right) \left( \partial_x z_j^{(P)} + \mu_i\left(u^{(P)}\right)\partial_t z_j^{(P)} \right)\\
  & +\sum_{j, k=1} ^n  B_{i j}\left(u^{(P)}\right) \frac{\partial \mu_i}{\partial u_k }\left(u^{(P)}\right) z_k^{(P)} z_j^{(P)}
\end{aligned}
$$
\begin{equation}\label{5.16}
\begin{aligned}
  & +\sum_{j, k=1}^n \frac{\partial B_{i j}}{\partial u_k}\left(u^{(P)}\right) z_k^{(P)}\left(\omega_j^{(P)}+\mu_i\left(u^{(P)}\right) z_j^{(P)}\right)\\
  & +\sum_{j \neq i} \tilde{g}_{i j} z_j^{(P)}
    +K \mu_i(0) z_i^{(P)}
    +\sum_{j=1}^n \frac{\partial \tilde{g}_i^{N L}}{\partial u_j}\left(u^{(P)}\right) z_j^{(P)},
\q \forall i=1,...,n.
\end{aligned}
\end{equation}
Noting that
$$
\begin{aligned}
& \sum_{j=1}^n B_{i j}(u)\left(\partial_x z_j+\mu_i(u) \partial_t z_j\right)-B_{i j}\left(u^{(P)}\right)\left(\partial_x z_j^{(P)}+\mu_i\left(u^{(P)}\right) \partial_t z_j^{(P)}\right) \\
= & \sum_{j=1}^n\left(\partial_x+\mu_i(u) \partial_t\right)\left(B_{i j}(u) z_j\right)-\left(\partial_x+\mu_i\left(u^{(P)}\right) \partial_t\right)\left(B_{i j}\left(u^{(P)}\right) z_j^{(P)}\right) \\
& -\sum_{j, j^{\prime}=1}^n z_j \frac{\partial B_{i j}}{\partial u_{j^{\prime}}}(u)\left(w_{j^{\prime}}+\mu_i(u) z_{j^{\prime}}\right)-z_j^{(P)} \frac{\partial B_{i j}}{\partial u_{j^{\prime}}}\left(u^{(P)}\right)\left(w_{j^{\prime}}^{(P)}+\mu_i\left(u^{(P)}\right) z_{j^{\prime}}^{(P)}\right),
\end{aligned}
$$
thus we have
\begin{equation}\label{5.17}
\begin{aligned}
   & \partial_x\left(z_i-z_i^{(P)}\right)+\mu_i(u) \partial_t\left(z_i-z_i^{(P)}\right)-\tilde{g}_{i i}\left( z_i - z_i^{(P)}\right) \\
=  &-\left(\mu_i(u)-\mu_i\left(u^{(P)}\right)\right) \partial_t z_i^{(P)}+\sum_{j=1}^n\left(\mu_i(u)-\mu_i\left(u^{(P)}\right)\right) B_{i j}\left(u^{(P)}\right) \partial_t z_j^{(P)} \\
& +\sum_{j=1}^n\left(\partial_x+\mu_i(u) \partial_t\right)\left(\left(B_{i j}(u)-B_{i j}\left(u^{(P)}\right)\right) z_j+B_{i j}\left(u^{(P)}\right)\left(z_j-z_j^{(P)}\right)\right) \\
& +\sum_{j, k=1}^n\left(\mu_i(u)-\mu_i\left(u^{(P)}\right)\right) \frac{\partial B_{i j}}{\partial u_{k}}\left(u^{(P)}\right) z_{k}^{(P)} z_j^{(P)} \\
& -\sum_{j, k =1}^n\left(\left(\frac{\partial B_{i j}}{\partial u_{k}}(u)-\frac{\partial B_{i j}}{\partial u_{k}}\left(u^{(P)}\right)\right) z_j w_{k}+\frac{\partial B_{i j}}{\partial u_{k}}\left(u^{(P)}\right)\left(z_j-z_j^{(P)}\right) w_{k}\right) \\
& -\sum_{j, k=1}^n \frac{\partial B_{i j}}{\partial u_{k}}\left(u^{(P)}\right) z_j^{(P)}\left(w_{k}-w_{k}^{(P)}\right)
  -\sum_{j=1}^n\left(\frac{\partial \mu_i}{\partial u_j}(u)-\frac{\partial \mu_i}{\partial u_j}\left(u^{(P)}\right)\right) z_j z_i \\
& -\sum_{j, k=1}^n\left(\frac{\partial B_{i j}}{\partial u_{k}}(u) \mu_i(u)-\frac{\partial B_{i j}}{\partial u_{k}}\left(u^{(P)}\right) \mu_i\left(u^{(P)}\right)\right) z_j z_{k}\\
& -\sum_{j, k=1}^n \frac{\partial B_{i j}}{\partial u_{k}}\left(u^{(P)}\right) \mu_i\left(u^{(P)}\right)\left(\left(z_j-z_j^{(P)}\right) z_{k}+z_j^{(P)}\left(z_{k}-z_{k}^{(P)}\right)\right) \\
&  -\sum_{j=1}^n \frac{\partial \mu_i}{\partial u_j}\left(u^{(P)}\right)\left(\left(z_j-z_j^{(P)}\right) z_i+z_j^{(P)}\left(z_i-z_i^{(P)}\right)\right)\\
& +\sum_{j, k=1}^n\left(B_{i j}(u) \frac{\partial \mu_i}{\partial u_{k}}(u)-B_{i j}\left(u^{(P)}\right) \frac{\partial \mu_i}{\partial u_{k}}\left(u^{(P)}\right)\right) z_{k} z_j\\
& +\sum_{j, k=1}^n B_{i j}\left(u^{(P)}\right) \frac{\partial \mu_i}{\partial u_{k}}\left(u^{(P)}\right)\left(\left(z_{k}-z_{k}^{(P)}\right) z_j+z_{k}^{(P)}\left(z_j-z_j^{(P)}\right)\right) \\
& +\sum_{j, k=1}^n\left(\frac{\partial B_{i j}}{\partial u_{k}}(u)-\frac{\partial B_{i j}}{\partial u_{k}}\left(u^{(P)}\right)\right) z_{k}\left(w_j+\mu_i(u) z_j\right) \\
& +\sum_{j, k=1}^n \frac{\partial B_{i j}}{\partial u_{k}}\left(u^{(P)}\right) \left(z_{k}-z_{k}^{(P)}\right)\left(w_j+\mu_i(u) z_j\right) \\
& +\sum_{j, k=1}^n \frac{\partial B_{i j}}{\partial u_{k}}\left(u^{(P)}\right) z_{k}^{(P)}\left(\left(w_j-w_j^{(P)}\right)+\mu_i(u)\left(z_j-z_j^{(P)}\right)\right) \\
& +\sum_{j, k=1}^n \frac{\partial B_{i j}}{\partial u_{k}}\left(u^{(P)}\right) z_{k}^{(P)}\left(\mu_i(u)-\mu_i\left(u^{(P)}\right)\right) z_j^{(P)}
  +\sum_{j \neq i} \tilde{g}_{i j} \left(z_j - z_j^{(P)} \right)\\
&  +K \mu_i(0)  \left(z_i - z_i^{(P)} \right)
  +\sum_{j=1}^n \frac{\partial \tilde{g}_i^{N L}}{\partial u_j}\left(u^{(P)}\right) \left(z_j - z_j^{(P)} \right)\\
& +\sum_{j=1}^n \left(\frac{\partial \tilde{g}_i^{N L}}{\partial u_j}\left(u\right) - \frac{\partial \tilde{g}_i^{N L}}{\partial u_j}\left(u^{(P)}\right)\right) z_j^{(P)},
  \q \forall i=1,...,n.
\end{aligned}
\end{equation}
For $r=1,...,m$ and $s=m+1,...,n$, we multiply $(\ref{5.17})$ with $W_r(x)$ and $W_s(x)$ respectively, and integrate along the corresponding characteristic curve $t=t_i\left(x ; \tilde{t}, \tilde{x}\right)$ defined by
\begin{equation}\label{5.18}
\left\{\begin{array}{l}
\displaystyle\frac{d}{d x} t_i\left(x ; \tilde{t}, \tilde{x}\right)
=\mu_i\left(u\left(t_i\left(x ; \tilde{t}, \tilde{x}\right), x\right)\right), \\\\
t_i\left(\tilde{x}; \tilde{t}, \tilde{x}\right)=\tilde{t}.
\end{array}\right.
\end{equation}
Noting the definition of $T_0$ in $(\ref{1.28})$, then for every point $\left(\tilde{t}, \tilde{x}\right) \in\left[(N+1) T_0, \tau\right]$, the characteristic curve $t=t_i\left(x ; \tilde{t}, \tilde{x}\right)$ will intersect the boundary in a time interval shorter than $T_0$,
i.e.,
$$
\begin{array}{ll}
t_r\left(L ;\tilde{t}, \tilde{x}\right) \in\left[ \tilde{t}-T_0,  \tilde{t}\right] \subseteq\left[N T_0, \tau\right],
 \quad \forall \tilde{t} \in\left[(N+1) T_0, \tau\right], \\
t_s\left(0 ;\tilde{t}, \tilde{x}\right) \in\left[ \tilde{t}-T_0,  \tilde{t}\right] \subseteq\left[N T_0, \tau\right],
 \quad \forall \tilde{t} \in\left[(N+1) T_0, \tau\right].
\end{array}
$$
Therefore, from $(\ref{2.3})$, $(\ref{2.8})$, $(\ref{5.1})$-$(\ref{5.4})$, and $(\ref{5.12})$-$(\ref{5.13})$ , for $r=1,...,m$, we have
\begin{equation}\label{5.19}
\begin{aligned}
& \left|z_r\left(\tilde{t}, \tilde{x}\right)-z_r^{(P)}\left(\tilde{t}, \tilde{x}\right)\right| \\
\leq &\frac{ \beta^{*} \widetilde{C}_S \varepsilon \beta^{N} }{W_r\left(\tilde{x}\right)}
      +\frac{\left(\sum\limits_{j\neq r}|\tilde{g}_{rj}|\right) \left(W_r(\tilde{x})-1\right) } {|\tilde{g}_{rr}| W_r\left(\tilde{x}\right)} \cdot \widetilde{C}_S \varepsilon \beta^N
       +\frac{C \left(\varepsilon + K\right)}{W_r\left(\tilde{x}\right)}\cdot \widetilde{C}_S \varepsilon \beta^N\\\\
     & +\max_{1=1,...,n}\sup_{u \in \mathcal{U}}|\nabla\mu_i(u)| M_R C_S \varepsilon \beta^{N}
       + C \varepsilon \cdot C_S \varepsilon \beta^{N}
        + C\varepsilon^2.
\end{aligned}
\end{equation}
Here we take the large enough constant $\widetilde{C}_S$ satisfying
\begin{equation}\label{5.20}
\widetilde{C}_S \geq \frac{\max\limits_{1=1,...,n}\sup\limits_{u \in \mathcal{U}}|\nabla\mu_i(u)| M_R C_S M_3}{1-\beta^{*}-(1-\beta) M_3} +C_S
\end{equation}
with $$(1-\beta)M_3 <1-\beta^{*}.$$
Then we have
\begin{equation}\label{5.20'}
 \left|z_r\left(\tilde{t}, \tilde{x}\right)-z_r^{(P)}\left(\tilde{t}, \tilde{x}\right)\right|
     \leq  \widetilde{C}_S \varepsilon \beta^{N+1},
     \quad \forall\left(\tilde{t}, \tilde{x}\right) \in\left[(N+1) T_0, \tau\right] \times[0, L].
\end{equation}
Similarly for $s=m+1,...,n$, we have
\begin{equation}\label{5.21}
 \left|z_s\left(\tilde{t}, \tilde{x}\right)-z_s^{(P)}\left(\tilde{t}, \tilde{x}\right)\right|
     \leq  \widetilde{C}_S \varepsilon \beta^{N+1},
     \quad \forall\left(\tilde{t}, \tilde{x}\right) \in\left[(N+1) T_0, \tau\right] \times[0, L].
\end{equation}
Combine $(\ref{2.8'})$  and $(\ref{3.11})$ to get
\begin{equation}\label{5.22}
\begin{aligned}
    w_i-w_i^{(P)}
= & \ \tilde{g}_{i i} \left( u_i -u_i^{(P)}\right)
    -\left(\mu_i(u)-\mu_i\left(u^{(P)}\right)\right) z_i-\mu_i\left(u^{(P)}\right)\left(z_i-z_i^{(P)}\right) \\
& +\sum_{j=1}^n\left(B_{i j}(u)-B_{i j}\left(u^{(P)}\right)\right) w_j+\sum_{j=1}^n B_{i j}\left(u^{(P)}\right)\left(w_j-w_j^{(P)}\right) \\
& +\sum_{j=1}^n\left(B_{i j}(u) \mu_i(u)-B_{i j}\left(u^{(P)}\right) \mu_i\left(u^{(P)}\right)\right) z_j \\
& +\sum_{j=1}^n B_{i j}\left(u^{(P)}\right) \mu_i\left(u^{(P)}\right)\left(z_j-z_j^{(P)}\right)\\
& +\sum_{j\neq i} \tilde{g}_{i j} \left( u_i -u_i^{(P)}\right)
  +K \mu_i(0) \left( u_i -u_i^{(P)}\right)\\
&   +\tilde{g}_i^{N L}\left(u\right)-\tilde{g}_i^{N L}\left(u^{(P)}\right).
\quad i=1, \ldots, n.
\end{aligned}
\end{equation}
From $(\ref{1.13})$, $(\ref{2.3})$,  $(\ref{5.1})$-$(\ref{5.4})$, and $(\ref{5.19})$, $(\ref{5.21})$, we have
\begin{equation}\label{5.23}
\begin{aligned}
\left\|w_i(t, \cdot)-w_i^{(P)}(t, \cdot)\right\|_{C^{\circ}}
 & \leq  C_S \beta^{N}\varepsilon \cdot C \varepsilon + \widetilde{C}_S \beta^{N} \varepsilon \cdot C \varepsilon
       + C \varepsilon^2 \\
& \leq \widetilde{C}_S \beta^{N+1}\varepsilon,
\quad \forall t \in\left[(N+1) T_0, \tau\right], \forall i=1, \ldots, n.
\end{aligned}
\end{equation}
Thus we complete the proof of Theorem 5.
\section{Acknowledgment}

Peng Qu is supported in part by NSFC Grants No. 12122104, 11831011, and Shanghai Science and Technology Programs 21ZR1406000, 21JC1400600, 19JC1420101. Huimin Yu is supported in part by NSFC Grant No. 12271310 and Natural Science Foundation of Shandong Province ZR2022MA088.

\end{document}